\documentclass[11pt, a4paper]{amsart}
\usepackage{amssymb, array, amsmath, amscd, pdfpages, enumerate, amsthm, setspace, mathtools}

\usepackage[utf8]{inputenc}
\usepackage{amssymb}
\usepackage{bm}
\usepackage{xcolor}
\usepackage{overpic}

\usepackage[colorlinks=true,allcolors=magenta]{hyperref}

\newcommand{\eq}[1]{\begin{align*}#1\end{align*}}
\newcommand{\ieq}[1]{$#1$}
\newcommand{\eqn}[1]{\begin{align}#1\end{align}}
\newcommand*{\abs} [1]{\lvert#1\rvert}
\newcommand*{\norm}[1]{\lVert#1\rVert}
\newcommand*{\set} [1]{\{#1\}}
\newcommand*{\setm}[2]{\{\,#1\mid#2\,\}}   
\newcommand*{\iprod}[2]{\langle#1,#2\rangle}    
\newcommand*{\Abs}[2][default]{\ifthenelse{\equal{#1}{default}}{\left\lvert#2\right\rvert}{\ldelim{#1}{\lvert}#2\rdelim{#1}{\rvert}}}
\newcommand*{\Norm}[2][default]{\ifthenelse{\equal{#1}{default}}{\left\lVert#2\right\rVert}{\ldelim{#1}{\lVert}#2\rdelim{#1}{\rVert}}}

\newcommand{\conj}[1]{\overline{#1}}

\newcommand{\ga}{\alpha}
\newcommand{\gb}{\beta}
\renewcommand{\gg}{\gamma}
\newcommand{\gd}{\delta}
\newcommand{\gl}{\lambda}
\newcommand{\gw}{\omega}
\newcommand{\gs}{\sigma}
\newcommand{\eps}{\varepsilon}
\newcommand*{\C}{{\mathbb{C}}}
\newcommand*{\R}{{\mathbb{R}}}
\newcommand*{\N}{{\mathbb{N}}}
\DeclareMathOperator{\re}{Re}

\DeclareMathOperator{\diag}{diag}

\DeclareMathOperator*{\essup}{ess\,sup}

\newcommand*{\List}[2][1]{\set{#1,\ldots,#2}}

\newcommand{\citel}[2]{\cite[#2]{#1}}
\newcommand{\keyterm}[1]{\emph{#1}}
\newcommand{\mc}[1]{\mathcal{#1}}

\newcommand*{\inv}{^{-1}}
\newcommand*{\Lp}[1][p]{L^{#1}}
\newcommand*{\Lploc}[1][p]{L_{\text{loc}}^{#1}}
\newcommand*{\lp}[1][p]{\ell^{#1}}
\newcommand*{\Lin}{{\mathcal{L}}}    
\newcommand{\Dom}{\mathcal{D}}
\newcommand{\ran}{\mathcal{R}}
\renewcommand{\ker}{\mathcal{N}}

\newcommand{\pmat}[1]{\begin{bmatrix}#1\end{bmatrix}}
\newcommand{\pmatsmall}[1]{\begin{bsmallmatrix}#1\end{bsmallmatrix}}

\newcommand{\potential}[1]{{\color{darkgreen} #1}}
\renewcommand{\potential}[1]{#1}

\newcommand{\cref}{c_{\mbox{\scriptsize\textit{ref}}}}

\newcommand{\yref}{y_{\mbox{\scriptsize\textit{ref}}}}

\newcommand{\yrefamplk}[1][k]{y_{\mbox{\scriptsize\textit{ref}}}^{#1}}
\newcommand{\wdistamplk}[1][k]{w_{\mbox{\scriptsize\textit{dist}}}^{#1}}

\newcommand{\wdist}{w_{\mbox{\scriptsize\textit{dist}}}}
\newcommand{\wdistk}[1][k]{d_{#1}}

\newcommand{\cextk}[1][k]{c^{#1}_e}

\newcommand{\FE}{\mu}
\newcommand{\Merr}{M_{err}}
\newcommand{\wext}{w_e}
\newcommand{\F}{\mathbb{F}}

\newcommand*{\ddb}[2][1]{\ifthenelse{\equal{#1}{1}}{\frac{d}{d#2}}{\frac{d^{#1}}{d#2^{#1}}}}
\newcommand*{\pd}[3][1]{\ifthenelse{\equal{#1}{1}}{\frac{\partial{#2}}{\partial{#3}}}{\frac{\partial^{#1}{#2}}{\partial#3^{#1}}}}
\newcommand*{\floor}[1]{\lfloor#1\rfloor}

\newcommand{\QUES}[1]{{\color{red} #1}}
\renewcommand{\QUES}[1]{}

\newcommand{\EAC}{\textup{EAC}}
\newcommand{\G}{\mathcal{G}}
\newcommand{\CL}{C_\Lambda}

\newcommand{\yaux}{y_{\mbox{\scriptsize\textit{aux}}}}
\newcommand{\yauxt}{\tilde{y}_{\mbox{\scriptsize\textit{aux}}}}
\newcommand{\ypert}{y_{\mbox{\scriptsize\textit{pert}}}}

\newcommand{\yauxrvec}{r}
\newcommand{\ytrans}{y_0}
\newcommand{\Kaux}{K_{\mbox{\scriptsize\textit{aux}}}}

\newcommand{\Ptot}{P_{\mbox{\scriptsize\textit{tot},$L$}}}

\newtheorem{theorem}{Theorem}[section]
\newtheorem{lemma}[theorem]{Lemma}

\theoremstyle{definition}
\newtheorem{definition}[theorem]{Definition}
\newtheorem{assumption}[theorem]{Assumption}

\newtheorem{remark}[theorem]{Remark}

\begin{document}

\title[Time-Varying Internal Models]{Time-Varying Internal Models and Regulation of Unknown Harmonic Signals for Regular Linear Systems}

\thispagestyle{plain}

\author[L. Paunonen and S. Afshar]{Lassi Paunonen and Sepideh Afshar}

\address[L. Paunonen]{Mathematics and Statistics, Tampere University, PO.\ Box 692, 33101 Tampere, Finland.}
\email{lassi.paunonen@tuni.fi}

\address[S. Afshar]{Gordon Center for Medical Imaging, Harvard Medical School, Massachusetts General Hospital, Boston, MA, USA.}
\email{safshar1@mgh.harvard.edu}

\thanks{The research was supported by the Academy of Finland Grants number 298182 and 310489 held by L. Paunonen.}

\begin{abstract}
We introduce general results on well-posedness and output regulation of regular linear systems with nonautonomous controllers. 
We present a generalization of the internal model principle for time-dependent controllers with asymptotically converging parameters.
This general result is utilised in controller design for output tracking and disturbance rejection of harmonic signals with unknown frequencies.
Our controller can be flexibly combined with different frequency estimation methods.
The results are illustrated in rejection of 
unknown harmonic disturbances 
for a one-dimensional boundary controlled heat equation.%
\end{abstract}


\maketitle

\section{Introduction}
\label{sec:introduction}

Asymptotic output tracking and disturbance rejection, jointly called output regulation, is an important control objective in many engineering applications.
For a given reference signal $\yref(t)$ and a class of external disturbance signals $\wdist(t)$
the output $y(t)$ of the controlled system is required to satisfy 
\eq{
  \norm{y(t)-\yref(t)}\to 0, \qquad \mbox{as} \quad t\to\infty.
}
This control problem 
has been studied extensively for both abstract infinite-dimensional systems~\cite{Sch83b,ByrLau00,HamPoh00,RebWei03,HamPoh10,NatGil14,Pau16a} and controlled partial differential equations (PDEs)~\cite{AnfStr17,DeuGab18,JinGuo19}.

In the classical output regulation problem the reference and disturbance signals are assumed to have the forms
\begin{subequations}
  \label{eq:yrefwdist}
  \eqn{
  \label{eq:yrefwdistyr}
    \hspace{-1ex}  \yref(t) &=
  \sum_{k=0}^q \yrefamplk \cos(\gw_k t+\theta_k) \\
  \label{eq:yrefwdistwd}
  \hspace{-1ex} \wdist(t) &=
  \sum_{k=0}^q \wdistamplk \cos(\gw_k t + \varphi_k),
  }
\end{subequations}
where the frequencies $0=\gw_0<\gw_1<\ldots<\gw_q$
 are assumed to be known 
and the
amplitudes $(\yrefamplk)_{k=0}^q$, $(\wdistamplk)_{k=0}^q$
 and phases $(\theta_k)_{k=0}^q$, $(\varphi_k)_{k=0}^q$  may be unknown.
In this paper we focus on a more challenging version of the control problem where also the frequencies $(\gw_k)_{k=1}^q$ are unknown. As our ultimate contribution we introduce 
 a controller design method for output regulation of signals~\eqref{eq:yrefwdist} with unknown frequencies, amplitudes and phases.
In particular, we introduce a dynamic error feedback controller with a \emph{time-varying internal model} based on estimates $(\hat{\gw}_k(t))_{k=1}^q$ of the 
frequencies 
 in~\eqref{eq:yrefwdist}.

We present our results for 
a \keyterm{regular linear system}~\cite{Wei94,Sta05book}%
\begin{subequations}%
\label{eq:RLSplant}
\eqn{
\dot{x}(t)&=Ax(t)+Bu(t) + B_d\wdist(t), \qquad x(0)=x_0\\
y(t)&=\CL x(t)+Du(t) + D_d \wdist(t),
\label{eq:RLSplantoutput}
}
\end{subequations}
on a Hilbert space $X$ (see Section~\ref{sec:prelimsMain} for detailed assumptions).
This class of systems covers a wide range of PDE models with boundary control and observation, e.g.,  one-dimensional convection-diffusion equations, wave equations, beam equations, as well as two-dimensional heat equations~\cite{ByrGil02}. 
Our main controller design method
for output regulation
leads to a \emph{nonautonomous} dynamic error feedback controller
\begin{subequations}
\label{eq:Contr}
\eqn{
  \hspace{-.2cm}\dot{z}(t)&= \mc{G}_1(t)z(t)+\mc{G}_2(t)(y(t)-\yref(t)), \quad z(0)=z_0\\
  \hspace{-.2cm}u(t)&= K(t)z(t)
}
\end{subequations}
on a Hilbert space $Z$.
Here 
 $\mc{G}_2(\cdot)\in \Lp[\infty](0,\infty;\Lin(\C^p,Z))$ and $K(\cdot)\in \Lp[\infty](0,\infty;$ $\Lin(Z,\C^m))$ and $\G_1(t)$ may contain an unbounded time-varying part (see Assumption~\ref{ass:ContrAss}).
The analysis of well-posedness of the closed-loop system consisting of~\eqref{eq:RLSplant} and~\eqref{eq:Contr} is highly nontrivial for unbounded operators $B$ and $C$.
As our first main result we prove that 
the closed-loop system has a well-defined mild state and output determined by bounded input/output maps. We achieve this result by expressing the time-varying closed-loop system as a \emph{nonautonomous output feedback} of an autonomous regular linear system and by employing the nonautonomous feedback theory developed by Schnaubelt in~\cite{Sch02}.
Besides output regulation, the well-posedness 
result in Section~\ref{sec:CLwellposedness} is 
also applicable in the study of other control problems with nonautonomous controllers.

In the second part of the paper we 
introduce general theory for output regulation
in the situation where the controller parameters $(\mc{G}_1(t),\mc{G}_2(t),K(t))$ --- especially the internal model --- 
 converge
 to a limit $(\mc{G}_1^\infty,\mc{G}_2^\infty,K^\infty)$ as $t\to \infty$. 
As our main result we show that if the 
autonomous ``limit controller'' $(\mc{G}_1^\infty,\mc{G}_2^\infty,K^\infty)$ contains an internal model~\cite{FraWon75a,PauPoh10} of the true frequencies
 $(\gw_k)_k$
 and the closed-loop system is exponentially stable, then the controller achieves output regulation. 
\potential{Moreover, we show that output tracking is achieved in an approximate sense
 even if the 
 limits of the frequencies in the internal model are only close to the frequencies 
$(\gw_k)_k$.
}

In Section~\ref{sec:ContrDesign} we introduce our controller for output regulation of signals~\eqref{eq:yrefwdist} with unknown frequencies.
We begin by introducing a general controller structure with a time-varying internal model and an observer part for closed-loop stabilization.
The controller also includes an auxiliary output $\yaux(t)$ for estimation of
 the frequencies $(\gw_k)_k$ in~\eqref{eq:yrefwdist}. 
One of the key features of our controller is that
$\yaux(t)$ is by design independent
of the time-varying parts of the controller, and therefore the convergence 
of the frequency estimates
 in the internal model can be completed \emph{separately} 
of the analysis of the closed-loop dynamics. 
In particular, our controller is not restricted to a single estimation method, but can instead 
be combined with 
any method which can identify $(\gw_k)_k$ based on the output $\yaux(t)$.
As the final part of the controller design we present an online tuning algorithm for the stabilization of the nonautonomous closed-loop system and for guaranteeing the output tracking.
In Section~\ref{sec:ContrRobAnalysis} we analyse the robustness properties of the controller with respect to perturbations in the system $(A,B,C,D)$. 
The detailed robustness properties depend on the chosen frequency  estimation method
due to the effect of the perturbations on $\yaux(t)$.
 Our main result shows that for sufficiently long update intervals
our controller achieves approximate output regulation despite small perturbations
 provided that the frequency estimates approximate the true frequencies 
$(\gw_k)_k$
 with sufficient accuracy.
We illustrate the controller design
 in Section~\ref{sec:simulation} in adaptive output regulation for a boundary controlled heat equation with uncertainty.

Output regulation of distributed parameter systems
with unknown frequencies 
$(\gw_k)_k$
 has been previously studied 
in~\cite{WanJi14,WanJi14b} where the system was assumed to be transformable into a canonical form.
Moreover, the output regulation of unknown harmonic signals has also been considered for particular PDEs, such as $2\times 2$ hyperbolic systems~\cite{AnfStr17},
  a Kirchoff plate~\cite{RitSch12}, and  a 1D boundary controlled heat equation~\cite{GuoZha22}.
This control problem has also been studied actively using adaptive internal models for finite-dimensional linear and nonlinear systems~\cite{SerIsi01,MarTom03,XuHua10,XuWan16,CarGal16}.
Our use of the auxiliary output $\yaux(t)$ is inspired by the ``residual generator'' 
in the internal model based controller in~\citel{CarGal16}{Sec.~4}.
A preliminary version of Theorem \ref{thm:ORPmain} was presented in~\cite{AfsPauACC19} for systems with bounded $B$ and $C$.

 \textbf{Notation.} If $X$ and $Y$ are Banach spaces, then the space of bounded linear operators $A:X \to Y$ is denoted by $\Lin(X,Y)$.
The domain, kernel, and range of $A : \Dom(A)\subset X\to Y$ are denoted by $\Dom(A)$, $\ker(A)$, and $\ran(A)$, respectively. The resolvent operator of $A: \Dom(A)\subset X\to X$ is defined as $R(\gl,A)=(\gl I-A)\inv$ for those $\gl\in \C$ for which the inverse is bounded.
The inner product on $X$ is denoted by $\iprod{\cdot}{\cdot}_X$.
By $\Lp[p](0,\tau;X)$ and $\Lp[\infty](0,\tau;X)$ we denote, respectively, the spaces of $p$-integrable and essentially bounded measurable functions $f:(0,\tau)\to X$.
For $f\in \Lp[\infty](0,\infty;X)$ we denote $\norm{f(t)}\to 0$ as $t\to \infty$ 
if 
 $\essup_{s\geq t}\norm{f(s)}\to 0$ as $t\to\infty$. 
If $A: \Dom(A)\subset X\to X$ generates a strongly continuous semigroup $T(t)$ on 
 $X$, we define
 $X_1=\Dom(A)$ equipped with the graph norm of $A$. Moreover, we define $X_{-1}$ as the completion of  $X$ with respect to the norm $\norm{x}_{-1}:=\norm{(\gl_0-A)\inv x}_X$ for a fixed $\gl_0\in\rho(A)$.
Then $A$ extends to an operator $ X\to X_{-1}$ (also denoted by $A$) and this extension generates a semigroup (also denoted by $T(t)$) on $X_{-1}$~\citel{TucWei09book}{Sec.~2.10}.

\section{Preliminaries and Standing Assumptions}
\label{sec:prelimsMain}

\subsection{Background on regular linear systems}

Let $X$, $U$, and $Y$ be Hilbert spaces and
 consider
\begin{subequations}
\label{eq:RLSprelim}
\eqn{
\dot{x}(t)&=Ax(t)+Bu(t),  \qquad x(0)=x_0\in X\\
y(t)&=\CL x(t)+Du(t) 
}
\end{subequations}
on $X$ with $B\in \Lin(U,X_{-1})$, $C\in \Lin(X_1,Y)$ and $D\in \Lin(U,Y)$.
Here 
 $C_\Lambda: \Dom(\CL)\subset X\to Y$ is the $\Lambda$-extension of $C$
such that
$\Dom(\CL) := \setm{x\in X}{\lim_{\gl\to\infty}\gl CR(\gl,A)x ~\mbox{exists}}$ 
and
\eq{
\CL x &= \lim_{\gl\to\infty}\gl CR(\gl,A)x , \qquad \forall x\in \Dom(\CL).
}
The operator $B$ is \emph{an admissible input operator} for the semigroup $T(t)$ if $\int_0^\tau T(t-s)Bu(s)ds\in X$ for all $u\in \Lp[2](0,\tau;U)$ and $\tau>0$~\citel{TucWei14}{Sec.~3}. Moreover, $C$ is \emph{an admissible output operator} for the semigroup $T(t)$ if there exist $\tau,\gg>0$ such that $\norm{CT(\cdot)x}_{\Lp[2](0,\tau)}\leq \gg \norm{x}$ for all $x\in \Dom(A)$~\citel{TucWei14}{Sec.~3}.

\begin{assumption}
\label{ass:RLS}
For some Hilbert spaces $X$, $U$ and $Y$ the operators $(A,B,C)$ have the following properties.
\begin{itemize}
\item[\textup{(1)}] $A: \Dom(A)\subset X\to X$ generates a semigroup $T(t)$ on $X$.
\item[\textup{(2)}] $B\in \Lin(U,X_{-1})$ and $C\in \Lin(X_1,Y)$ are admissible.
\item[\textup{(3)}] $\ran(R(\gl_0,A)B)\subset \Dom(\CL)$ for some $\gl_0\in \rho(A)$.
\item[\textup{(4)}]  $\sup_{\re\gl\geq \gb} \norm{\CL R(\gl,A)B}<\infty$ for some $\gb>0$.
\end{itemize}
\end{assumption}

If  $(A,B,C)$ satisfy Assumption~\ref{ass:RLS} and $D\in \Lin(U,Y)$, then~\eqref{eq:RLSprelim} is a \emph{regular linear system} by~\citel{TucWei14}{Thm.~5.6} and its transfer function is given by $P(\gl)=\CL R(\gl,A)B+D$. 
In this situation we write
 ``$(A,B,C,D)$ is a regular linear system''.

\subsection{Assumptions on the system and the controller}

We consider a regular linear system of the form~\eqref{eq:RLSplant}
on a Hilbert space $X$, where
$x(t)\in X$, $u(t)\in\C^m$, $y(t)\in\C^p$, and $\wdist(t)\in \C^{n_d}$
are the system's state, input, output, and external disturbance, respectively.
In particular, the number of outputs of the system is $p\in\N$.
The operator $A: D(A)\subset X\rightarrow X$ is assumed to generate a strongly continuous semigroup $T(t)$ on $X$ and
$B\in \Lin(\C^m,X_{-1})$, $B_d\in \Lin(\C^{n_d},X)$, $C\in \Lin(X_1,\C^p)$, $D_d\in \C^{p\times n_d}$.
The operators $B_d\in\Lin(\C^{n_d},X)$ and $D_d\in \C^{p\times n_d}$ are allowed to be unknown.
We assume that $(A,B,C,D)$ satisfy Assumption~\ref{ass:RLS}. Since $B_d\in \Lin(\C^{n_d},X)$, also $(A,[B,B_d],C,[D,D_d])$ is a regular linear system.
The transfer function of~\eqref{eq:RLSplant}
(from $u$ to $y$) is denoted by $P(\gl)=\CL R(\gl,A)B+D$ for $\gl\in\rho(A)$.

  We make the following assumptions on the parameters of the dynamic error feedback controller~\eqref{eq:Contr}.

\begin{assumption}
\label{ass:ContrAss}
For almost every $t\geq 0$ we have
\eq{
 \mc{G}_1(t)z&=\mc{G}_1^\infty z + \mc{G}_{11}^\infty \Delta_{\mc{G}_1}(t)z\\
\Dom(\mc{G}_1(t))&= \setm{z\in Z}{\mc{G}_1^\infty z + \mc{G}_{11}^\infty \Delta_{\mc{G}_1}(t)z\in Z},
}
 where $\mc{G}_1^\infty : \Dom(\mc{G}_1^\infty)\subset Z\to Z$ generates a strongly continuous semigroup
on $Z$, $\mc{G}_{11}^\infty\in \Lin(U_c,Z_{-1})$ for some Hilbert space $U_c$ is an admissible input operator for this semigroup,
 and $\Delta_{\mc{G}_1}(\cdot)\in \Lp[\infty](0,\infty;\Lin(Z,U_c))$.
Moreover, $\mc{G}_2(\cdot)\in \Lp[\infty](0,\infty;\Lin(\C^p,Z))$ and $K(\cdot)\in \Lp[\infty](0,\infty;\Lin(Z,\C^m))$.
\end{assumption}

\potential{In Section~\ref{sec:ContrDesign} the controller will contain
 additional dynamics for estimation of  $(\gw_k)_k$ in~\eqref{eq:yrefwdist}, but
in our control scheme
the convergence of the frequency estimates
is
 analysed 
separately.%
}

We can formally express the closed-loop system 
of~\eqref{eq:RLSplant} and~\eqref{eq:Contr}
with state $x_e(t)=[x(t),z(t)]^T\in X_e:= X\times Z$ and input $w_e(t)=[\wdist(t),\yref(t)]^T\in \C^{n_d+p}$
as
\begin{subequations}
\label{eq:CLS}
\eqn{
\label{eq:CLS1}
  \dot{x}_e(t)&= A_e(t)x_e(t)+B_e(t)\wext(t), \quad x_e(0)=x_{e0}\\
  e(t)&= C_e(t)x_e(t) + D_e\wext(t)
}
\end{subequations}
where $e(t)=y(t)-\yref(t)$, $D_e=[D_d,\;-I]\in \C^{p\times (n_d+p)}$,  
\eq{
  A_e(t)&=\pmat{A&BK(t)\\\mc{G}_2(t)\CL&\mc{G}_1(t)+\mc{G}_2(t)D K(t)}
\\
   B_e(t) &= \pmat{B_d&0\\\mc{G}_2(t)D_d&\;-\mc{G}_2(t)}, \,
  C_e(t)=\pmat{\CL, DK(t)},
}
and $\Dom(A_e(t))=\setm{[x,z]^T\in \Dom(\CL)\times \Dom(\mc{G}_1(t))}{Ax+BK(t)z\in X}$.
The existence of well-defined mild state $x_e(t)$ and output $e(t)$ of~\eqref{eq:CLS} are proved in Section~\ref{sec:CLwellposedness}.

\section{Well-Posedness of the Closed-Loop System}
\label{sec:CLwellposedness}

In this section we will prove that the closed-loop system is well-posed in the sense that
for the initial state $x_{e0}=[x_0,z_0]^T\in X_e$ and for $\wext\in \Lploc[2](0,\infty;\C^{n_d+p})$ the equations~\eqref{eq:CLS}
have a well-defined mild state and output%
\begin{subequations}
\label{eq:CLSstateoutput}
\eqn{
\label{eq:CLSstateoutputstate}
x_e(t) &= U_e(t,0)x_{e0} + \Phi_e^{t,0}\wext\\
e(t) &= (\Psi_e^0 x_{e0})(t) + (\F_e^0 \wext)(t),
}
\end{subequations}
where $U_e(t,s)$ is a strongly continuous evolution family~\citel{Sch02}{Def.~2.1} and $(U_e,\Phi_e,\Psi_e,\F_e)$ is a \emph{well-posed nonautonomous linear system}~\citel{Sch02}{Def.~3.6}.
We will show that $U_e(t,s)$ is related to $(A_e(t))_{t\geq 0}$ through a natural perturbation formula and that the mappings%
\begin{subequations}
\label{eq:CLMapsDef}
\eqn{
\label{eq:CLMapsDefPhi}
\Phi_e^{t,s} \wext &= \int_s^t U_e(t,r)B_e(r) \wext(r)dr\\
(\Psi_e^s x)(t) &= C_e(t) U_e(t,s)x\\
\label{eq:CLMapsDefF}
\hspace{-1ex}(\F_e^s \wext)(t)  &= C_e(t) \hspace{-.4ex}\int_s^t \hspace{-.8ex} U_e(t,r)B_e(r) \wext(r)dr + D_e\wext(t)
}
\end{subequations}
are well-defined for all $s\geq 0$ and a.e. $t\geq s$ and have appropriate boundedness properties
 (see Theorem~\ref{thm:CLWellposedness}).

\begin{figure}[h!]
\begin{center}
\includegraphics[width=0.55\linewidth]{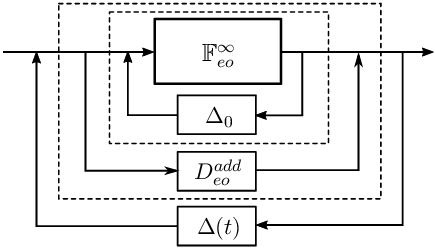}
\caption{The nonautonomous feedback structure.}
\label{fig:NRLSfeedback}
\end{center}
\end{figure}

We prove the closed-loop well-posedness using the nonautonomous feedback theory 
in~\citel{Sch02}{Sec.~4}.
More precisely, we will express the system~\eqref{eq:CLS} as a part of a system obtained from an autonomous regular linear system $(A_{eo}^\infty,B_{eo}^\infty,C_{eo}^\infty,D_{eo}^\infty)$ under
 a combination of 
(i) autonomous output feedback with feedback operator $\Delta_0$,
  (ii) parallel interconnection with a feedthrough operator $D_e^{add}$
and (iii) nonautonomous feedback with feedback operators $\Delta(t)$ (see Figure~\ref{fig:NRLSfeedback}).

In view of Section~\ref{sec:convfreq}, we let $\mc{G}_2^\infty\in \Lin(\C^p,Z)$ and $K^\infty\in \Lin(Z,\C^m)$ and define $\Delta_{\mc{G}_2}(\cdot) = \mc{G}_2(\cdot)-\mc{G}_2^\infty \in \Lp[\infty](0,\infty;\Lin(\C^p,Z))$ and $\Delta_{K}(\cdot)=K(\cdot)-K^\infty \in \Lp[\infty](0,\infty;\Lin(Z,\C^m))$.
We denote $U_{ee}=U\times U_d \times Y\times U_c\times Z$ and $Y_{ee} = Y\times U\times Z$ and 
 define $\Dom(A_{eo}^\infty)=\Dom(A)\times \Dom(\mc{G}_1^\infty)$,
$\Dom(C_{eo}^\infty)=\Dom(\CL)\times Z$,
\eq{
A_{eo}^\infty &= \pmat{A&0\\0&\mc{G}_1^\infty }, \quad
B_{eo}^\infty = \pmat{B&B_d&0&0&0\\0&0&-\mc{G}_2^\infty&\mc{G}_{11}^\infty&I}
\\
C_{eo}^\infty &= \pmat{\CL&0\\0&K^\infty\\0&I},
 \qquad
 D_{eo}^\infty =\pmat{D&D_d&0&0&0\\0&0&0&0&0\\0&0&0&0&0}.
}
 Our assumptions on $(A,[B,B_d],C,[D,D_d])$ and $\mc{G}_{11}^\infty$ imply that
 $(A_{eo}^\infty, B_{eo}^\infty,$ $ C_{eo}^\infty, D_{eo}^\infty)$ 
is a regular linear system with input space $U_{ee}$ and output space $Y_{ee}$. 
The operator $C_{eo}^\infty: \Dom(C_{eo}^\infty)\subset X_e\to Y_{ee}$ coincides with the $\Lambda$-extension of $C_{eo}^\infty: \Dom(A_{eo}^\infty)\subset X_e\to X_e$.
We define $\Delta_0\in \Lin(Y_{ee},U_{ee})$ and $D_{eo}^{add}\in \Lin(Y_{ee},Y_{ee})$ by
\eqn{
\label{eq:Delta0Deo}
\Delta_0 = \pmat{0&I&0\\0&0&0\\-I&0&0\\0&0&0\\0&0&0}, \quad
D_{eo}^{add} =\pmat{0&0&-I&0&0\\I&0&0&0&0\\0&0&0&0&0}.
}
Since $I-D_{eo}^\infty \Delta_0\in \Lin(Y_{ee})$ is boundedly invertible, $\Delta_0$ is an admissible output feedback operator for 
 $(A_{eo}^\infty, B_{eo}^\infty, C_{eo}^\infty, D_{eo}^\infty)$. The results in~\cite{Wei94} and~\citel{Sta05book}{Sec.~7.5} imply that first applying output feedback with operator $\Delta_0$ and subsequently adding a parallel connection
with the (constant) transfer function $D_{eo}^{add}$ produces a regular linear system $(A_e^\infty,B_{ee}^\infty,C_{ee}^\infty,D_{ee}^\infty)$ with 
\begin{subequations}
\label{eq:CLlimitoperators}
\eqn{
\label{eq:CLlimitoperatorsAe}
A_e^\infty &= A_{eo}^\infty + B_{eo}^\infty \Delta_0(I-D_{eo}^\infty\Delta_0)\inv C_{eo}^\infty\\
B_{ee}^\infty &=  B_{eo}^\infty (I-\Delta_0 D_{eo}^\infty)\inv\\
C_{ee}^\infty &=   (I-D_{eo}^\infty\Delta_0)\inv C_{eo}^\infty\\
D_{ee}^\infty &=   (I-D_{eo}^\infty\Delta_0)\inv  D_{eo}^\infty+D_{eo}^{add}.
}
\end{subequations}
where $\Dom(A_e^\infty) = \setm{x_e\in \Dom(\CL)\times Z}{A_e^\infty x_e\in X_e}$ 
and $\Dom(C_{ee}^\infty)=\Dom(C_{eo}^\infty)=\Dom(\CL)\times Z$.
We denote
 by $T_e(t)$ the strongly continuous semigroup generated by $A_e^\infty$  and  by $\F_{ee}^\infty$ 
the extended input-output map of $(A_e^\infty,B_{ee}^\infty,C_{ee}^\infty,D_{ee}^\infty)$.

We define $\Delta(\cdot)\in \Lp[\infty](0,\infty;\Lin(Y_{ee},U_{ee}))$ and $P_{in}\in \Lin(\C^{n_d+p},U_{ee})$ by
\eqn{
\label{eq:Deltat}
\Delta (t)=
\pmat{0&0&\Delta_K(t)\\0&0&0\\0&0&0\\0&0&\Delta_{\mc{G}_1}(t)\\\Delta_{\mc{G}_2}(t)&0&0}, 
\quad P_{in}=\pmat{0&0\\I&0\\0&I\\0&0\\0&0} 
}
and define $P_{out}=\bigl[I,0,0\bigr]\in \Lin(Y_{ee},Y)$.
A direct computation shows that for a.e. $t\geq 0$
the operator $I-D_{ee}^\infty\Delta(t)\in \Lin(Y_{ee})$ is boundedly invertible and
\eq{
 A_e^\infty + B_{ee}^\infty \Delta(t)(I-D_{ee}^\infty\Delta(t))\inv C_{ee}^\infty 
&= A_e(t),
}
where $A_e(t)$ is as in~\eqref{eq:CLS}.
This identity confirms that $A_e(t)$ are (at this stage formally) associated to the system obtained from $(A_e^\infty,B_{ee}^\infty,C_{ee}^\infty,D_{ee}^\infty)$ with the nonautonomous feedback $\Delta(t)$.
The following lemma shows that $\Delta(t)$ is an admissible feedback for  $(A_e^\infty,B_{ee}^\infty,C_{ee}^\infty,D_{ee}^\infty)$ in the sense of~\citel{Sch02}{Def.~4.1}.

\begin{lemma}
\label{lem:NAfeedbackadmissible}
Let $\F_{ee}^\infty$ be the extended input-output map of $(A_e^\infty,B_{ee}^\infty,C_{ee}^\infty,D_{ee}^\infty)$ and 
let $\Delta(\cdot)$ be as in~\eqref{eq:Deltat}.
Then
for every $t_0>0$ the operators
$I-\F_{ee}^\infty \Delta(\cdot)\in \Lin(\Lp[2](s,s+t_0;Y_{ee}))$ 
for $s\geq 0$
have uniformly bounded inverses.
\end{lemma}

\begin{proof}
The input-output map $\F_{eo}$ of $(A_{eo}^\infty, B_{eo}^\infty, C_{eo}^\infty, D_{eo}^\infty)$ 
can be partitioned as
\eq{
\F_{eo}
= \pmat{\F&\F_d&0&0&0\\
0&0&\F_c&\F_{c2}&\F_{c3}\\
0&0&\F_{c4}&\F_{c5}&\F_{c6} }.
}
We have $\F_{ee}^\infty
 = (I-\F_{eo} \Delta_0)\inv \F_{eo} +D_{eo}^{add} $, which implies
\eq{
\MoveEqLeft[1] I- \F_{ee}^\infty \Delta (\cdot)
=
\pmat{I&0&-S_1\F(\Delta_K(\cdot)+\F_{c2}\Delta_{\mc{G}_1}(\cdot))\\
0&I& -S_2(\Delta_K(\cdot)+\F_{c2}\Delta_{\mc{G}_1}(\cdot))\\
0&0&I}\\
&\quad -\pmat{\F S_2\F_{c3}\Delta_{\mc{G}_2}(\cdot)&0&0\\
S_2\F_{c3}\Delta_{\mc{G}_2}(\cdot)&0&  0\\
S_3
\Delta_{\mc{G}_2}(\cdot)&0& 
\F_{c4} S_1\F\Delta_K(\cdot)+S_4\Delta_{\mc{G}_1}(\cdot)
}
}
with $S_1=(I+\F\F_c)\inv$, $S_2=(I+\F_c\F)\inv$, 
$S_3=\F_{c6}-\F_{c4}\F S_2\F_{c3}$, and 
$S_4=\F_{c5}-\F_{c4}\F S_2\F_{c2}$.
The first term 
of $I- \F_{ee}^\infty \Delta (\cdot)$
is invertible on $\Lp[2](s,s+t_0;Y_{ee})$ and the inverse is uniformly bounded with respect to $s\geq 0$ and $0<t_0\leq 1$.
Since $K^\infty$ and $I$
 in $C_{eo}^\infty$ 
are bounded operators, it is easy to verify that the  
restrictions of the input-output maps
 $\F_{c3}$, $\F_{c4}$, $\F_{c5}$, and $\F_{c6}$ to the time-interval $[0,t_0]$ satisfy
$\norm{\F_{c3}\vert_{[0,t_0]}}\to 0$, 
$\norm{\F_{c4}\vert_{[0,t_0]}}\to 0$, 
$\norm{\F_{c5}\vert_{[0,t_0]}}\to 0$,  and
$\norm{\F_{c6}\vert_{[0,t_0]}}\to 0$ as $t_0\to 0$.
Since $\Delta_{\mc{G}_1}$, $\Delta_{\mc{G}_2}$ and $\Delta_K$ are essentially bounded
the $\Lin(\Lp[2](s,s+t_0;Y_{ee}))$-operator
norm of the second term of $I-\F_{ee}^\infty\Delta(\cdot)$ converges to zero as $t_0\to 0$ uniformly with respect to $s\geq 0$. 
Because of this, 
for a sufficiently small $t_0>0$ the operators 
$I-\F_{ee}^\infty \Delta(\cdot)\in \Lin(\Lp[2](s,s+t_0;Y_{ee}))$ for $s\geq 0$ have uniformly bounded inverses $\Lp[2](s,s+t_0;Y_{ee})$. By~\citel{Sch02}{Lem.~4.2} the same is then true for all $t_0>0$.
\end{proof}

We define 
 $\Dom(C_{ee}(t)) = \Dom(C_{ee}^\infty)=\Dom(\CL)\times Z$ for a.e. $t\geq 0$ and
\eqn{
\label{eq:BeeCeeDef}
B_{ee}(\cdot) = B_{ee}^\infty Q_1(\cdot) \quad \mbox{and} \quad
C_{ee}(\cdot) = Q_2(\cdot)C_{ee}^\infty
}
with
  $Q_1(\cdot) = (I-\Delta(\cdot)D_{ee}^\infty)\inv$ and $Q_2(\cdot) = (I-D_{ee}^\infty\Delta(\cdot))\inv$.
The following theorem shows that the closed-loop system~\eqref{eq:CLS} has a well-defined strongly continuous evolution family $U_e(t,s)$ and an input map $\Phi_e^{t,s}$ defined in~\eqref{eq:CLMapsDefPhi}. 
Moreover, since
a direct computation shows that 
$B_{ee}(t)P_{in} = B_e(t)$ and $P_{out}C_{ee}(t)=C_e(t)$ for a.e. $t\geq 0$, 
 the mappings defined in~\eqref{eq:CLMapsDef} and in the following theorem satisfy $\Psi_e^s=P_{out}\Psi_{ee}^s$ and $\F_e^s = P_{out}\F_{ee}^s$. Because of this, Theorem~\ref{thm:CLWellposedness} implies that for $x_0\in X$, $z_0\in Z$, $\wdist(t)$ and $\yref(t)$ the closed-loop system~\eqref{eq:CLS} has a well-defined mild state $x_e(t)$ and output $e(\cdot)\in \Lploc[2](0,\infty;\C^p)$ determined by~\eqref{eq:CLSstateoutput}.
The integral equation associates $(A_e(t))_{t\geq 0}$ to the evolution family $U_e(t,s)$.
 
\begin{theorem}
\label{thm:CLWellposedness}
Let Assumption~\textup{\ref{ass:ContrAss}} hold and
let  $\Delta(\cdot)$ be as in~\eqref{eq:Deltat}. 
There exists a strongly continuous evolution family $U_e(t,s)$ such that for all $x\in X_e$ 
and $s\geq 0$
 we have 
$U_e(r,s)x\in \Dom(C_{ee}^\infty)=\Dom(\CL)\times Z$
 for a.e. $r\geq s$, 
$\norm{C_{ee}^\infty U_e(s+\cdot,s)x}_{\Lp[2](s,s+t_0)}\leq \gg(t_0)\norm{x}$ for every $t_0>0$ and some $\gg(t_0)>0$ (depending only on $t_0>0$), and
\eq{
\MoveEqLeft[0.5] U_e(t,s)x = T_e(t-s)x
 + \hspace{-.4ex}\int_s^t \hspace{-.7ex} T_e(t-r)B_{ee}^\infty \Delta(r)(I-D_{ee}^\infty \Delta(r))\inv C_{ee}^\infty U_e(r,s)xds
}
for all $t\geq s$.
If $\Phi_e^{t,s}$ is defined as in~\eqref{eq:CLMapsDefPhi} and if we define
$(\Psi_{ee}^s x)(t) = C_{ee}(t) U_e(t,s)x$ and
\eq{
(\F_{ee}^s \wext)(t)  &= C_{ee}(t) \hspace{-.5ex}\int_s^t \hspace{-1.1ex} U_e(t,r)B_e(r) \wext(r)dr + D_{ee}^\infty P_{in}\wext(t),
}
then $(U_e,\Phi_e^{t,s},\Psi_{ee}^s,\F_{ee}^s)_{t\geq s\geq 0}$ is a well-posed nonautonomous system in the sense of~\textup{\citel{Sch02}{Def.~3.6}}.
 In particular,
 for all $s\geq 0$ and $t_0>0$ 
we have $\Phi_e^{t,s}\wext\in \Dom(C_{ee}(t)) = \Dom(\CL)\times Z$ a.e. $t\geq s$ and $\wext\in \Lploc[2](0,\infty;\C^{n_d+p})$,
\eq{
\Phi_e^{t,s}&\in \Lin(\Lp[2](s,t;\C^{n_d+p}),X_e), \qquad 0\leq t-s\leq t_0\\
\Psi_{ee}^s &\in \Lin(X_e,\Lp[2](s,s+t_0;Y)),\\
\F_{ee}^s&\in \Lin(\Lp[2](s,s+t_0;\C^{n_d+p}),\Lp[2](s,s+t_0;Y))
}
 with bounds independent of $s\geq 0$.
Finally, we have
\eq{
\F_{ee}^s  = (I-\F_{ee}^\infty\Delta(\cdot))\inv \F_{ee}^\infty P_{in}
}
where $\F_{ee}^\infty$ is the input-output map of $(A_e^\infty,B_{ee}^\infty,C_{ee}^\infty,D_{ee}^\infty)$.
\end{theorem}

\begin{proof}
Due to the regularity of $(A_e^\infty,B_{ee}^\infty,C_{ee}^\infty,D_{ee}^\infty)$, $B_{ee}(\cdot)$ and $C_{ee}(\cdot)$ defined in~\eqref{eq:BeeCeeDef} are ``admissible input and output operators for the evolution family $(T_e(t-s))_{t\geq s\geq 0}$'' in the sense of~\citel{Sch02}{Def.~3.3 and~2.4}, where  $T_e(t)$ is the semigroup generated by $A_e^\infty$.
For $t\geq s\geq 0$ and $u\in \Lp[2](s,t;U_{ee})$ we define the mapping $\overline{\mathbb{K}}_s$ as in~\citel{Sch02}{Def.~3.3} by
\eq{
(\overline{\mathbb{K}}_s B_{ee}(\cdot)u)(t) := \int_s^t T_e(t-r)B_{ee}(r)u(r)dr.
}
The regularity of $(A_e^\infty,B_{ee}^\infty,C_{ee}^\infty,D_{ee}^\infty)$ implies that 
$(\overline{\mathbb{K}}_s B_{ee}(\cdot)u)(t)\in \Dom(C_{ee}^\infty)=\Dom(C_{ee}(t))$ for a.e. $t\geq s$ and $C_{ee}(\cdot)\overline{\mathbb{K}}_s B_{ee}(\cdot)u\in \Lploc[2](s,\infty;Y_{ee})$.
Moreover, if we define $\Phi_0^{t,s}$, $\Psi_0^s$ and $\F_0^s$ by
\eq{
\Phi_{0}^{t,s}u = (\overline{\mathbb{K}}_s B_{ee}(\cdot)u)(t), \quad
(\Psi_{0}^s x)(t) = C_{ee}(t)T_e(t-s)x,
}
and $\F_0^s u = C_{ee}(\cdot)\overline{\mathbb{K}}_s B_{ee}(\cdot)u$
for $u\in \Lp[2](s,t;U_{ee})$ and $x\in X_e$, then 
$(T_e,\Phi_0^{t,s},$ $\Psi_0^{s},\F_0^s)_{t\geq s\geq 0}$ is a well-posed nonautonomous system by~\citel{Sch02}{Lem.~3.9}. 

We will now show that $(I-\Delta(\cdot)D_{ee}^\infty)\Delta(\cdot)$ is an admissible feedback~\citel{Sch02}{Def.~4.1} for
$(T_e,\Phi_0^{t,s},\Psi_0^{s},\F_0^s)_{t\geq s\geq 0}$. Due to the definitions, we have $\F_0^s=Q_2(\cdot)[\F_{ee}^\infty - D_{ee}^\infty]Q_1(\cdot)$ for all $s\geq 0$
where  $Q_1(\cdot) = (I-\Delta(\cdot)D_{ee}^\infty)\inv$ and $Q_2(\cdot) = (I-D_{ee}^\infty\Delta(\cdot))\inv$.
 Thus
\eq{
\MoveEqLeft[5]
I-\F_0^s (I-\Delta(\cdot)D_{ee}^\infty)\Delta(\cdot)
= I -Q_2(\cdot) (\F_{ee}^\infty - D_{ee}^\infty)\Delta(\cdot)\\
&= I-Q_2(\cdot)\F_{ee}^\infty \Delta(\cdot) + (I-D_{ee}^\infty\Delta(\cdot))\inv D_{ee}^\infty \Delta(\cdot)\\
&= Q_2(\cdot) (I-\F_{ee}^\infty \Delta(\cdot)).
}
Since $\Delta(\cdot)$ is an admissible feedback for $\F_{ee}^\infty$ by Lemma~\ref{lem:NAfeedbackadmissible} and since $Q_2(\cdot) $ have uniformly bounded inverses on $\Lp[\infty](0,\infty;\Lin(Y_{ee}))$, 
we have that $(I-\Delta(\cdot)D_{ee}^\infty)\Delta(\cdot)\in \Lp[\infty](0,\infty;\Lin(Y_{ee},U_{ee}))$ is an admissible feedback for $(T_e,\Phi_0^{t,s},$ $\Psi_0^{s},\F_0^s)_{t\geq s\geq 0}$.

By construction,  
the system $(T_e,\Phi_0^{t,s},\Psi_0^{s},\F_0^s)_{t\geq s\geq 0}$ has the properties in the first part of~\citel{Sch02}{Thm.~3.11}, namely,
that 
$\Phi_{0}^{t,s}u\in \Dom(C_{ee}(t))$ for a.e. $t\geq s$ and $t\mapsto C_{ee}(t)\Phi_{0}^{t,s}u\in \Lploc[2](s,\infty;Y_{ee})$ for all 
$s\geq 0$ and $u\in \Lploc[2](s,\infty;U_{ee})$. 
Because of this, the proof of~\citel{Sch02}{Thm.~4.4(a)} shows that there exists 
a strongly continuous evolution family $U_e(t,s)$ which satisfies the integral equation in the claim, and $C_{ee}^\infty U_e(\cdot,s)x\in \Lploc[2](s,\infty;Y_{ee})$ for all $x\in X_e$ and $s\geq 0$.
The proof of~\citel{Sch02}{Thm.~4.4(a)} also shows that $C_{ee}(t) = (I-D_{ee}^\infty\Delta(t))\inv C_{ee}^\infty$, a.e. $t\geq 0$, are admissible observation operators for $U_e(t,s)$.
If 
$\Psi_{ee}^s$ is defined as in the claim, then 
by~\citel{Sch02}{Lem.~2.5} the pair
 $(U_e,\Psi_{ee}^s)$
is a ``nonautonomous observation system'' in the sense of~\citel{Sch02}{Def.~2.2}.
Since $B_e(t)\in \Lin(\C^{n_d+p},X_e)$ for a.e. $t\geq 0$,  $\Phi_e^{t,s}$ can be defined as in~\eqref{eq:CLMapsDefPhi}, and $(U_e,\Phi_e^{t,s})$ is a ``nonautonomous control system'' in the sense of~\citel{Sch02}{Def.~3.1}. Since $B_e(\cdot)\in \Lp[\infty](0,\infty;\Lin(\C^{n_d+p},X_e))$, the properties of $(U_e,\Psi_e^s)$ and~\citel{Sch02}{Prop.~2.11} also imply that
there exist $\kappa,t_1>0$ such that
 for all $s\geq 0$ and $\wext\in \Lploc[2](0,\infty;\C^{n_d+p})$ we have $\Phi_e^{t,s}\wext\in \Dom(C_{ee}(t))$ for a.e. $t\geq s$ and
\eq{
\norm{C_{ee}(\cdot)\Phi_e^{\cdot,s}\wext}_{\Lp[2](s,s+t_1)}\leq \kappa \norm{\wext}_{\Lp[2](s,s+t_1)}.
}
Thus~\citel{Sch02}{Lem.~3.9} implies that
 $(U_e,\Phi_e^{t,s},\Psi_{ee}^s,\F_{ee,0}^s)_{t\geq s\geq 0}$ with $\F_{ee,0}^s:=C_e(\cdot)\Phi_e^{\cdot,s}$ is a well-posed nonautonomous system, and since $\F_{ee}^s\wext = \F_{ee,0}^s \wext+D_{ee}^\infty P_{in}\wext$, the same is finally true also for $(U_e,\Phi_e^{t,s},\Psi_{ee}^s,\F_{ee}^s)_{t\geq s\geq 0}$. In particular, $\Phi_e^{t,s}$, $\Psi_{ee}^s$, $\F_{ee}^s$ have the boundedness properties in the claim.

Finally, we will show that 
$\F_{ee}^s  = (I-\F_{ee}^\infty\Delta(\cdot))\inv \F_{ee}^\infty P_{in}$.
Let $\wext\in \Lploc[2](0,\infty;\C^{n_d+p})$.
The evolution family $U_e(t,s)$ is associated to  $(U_e,\Phi_e^{t,s},$ $\Psi_{ee}^s,\F_{ee,0}^s)_{t\geq s\geq 0}$ which is obtained from 
$(T_e,\Phi_0^{t,s},\Psi_0^{s},\F_0^s)_{t\geq s\geq 0}$
with output feedback $(I-\Delta(\cdot)D_{ee}^\infty)\Delta(\cdot)$.
Applying the identity (4.13) in the proof of~\citel{Sch02}{Thm.~4.4(b)}\footnote{The identity (4.13) does not require ``absolute regularity''
and it extends to $\Lploc[2](0,\infty;X_e)$ since 
$(U_e,\Psi_{ee}^s)$
is a nonautonomous observation system.} to these two systems and 
 $f=B_e(\cdot)\wext\in \Lploc[2](0,\infty;X_e)$ 
and using $B_e(\cdot)=B_{ee}(\cdot)P_{in}$ 
we get%
\eq{
\F_{ee}^s \wext \hspace{-.4ex} -\hspace{-.4ex} D_{ee}^\infty P_{in} \wext
\hspace{-.25ex}&=\hspace{-.15ex}\F_{ee,0}^s\wext
\hspace{-.25ex}= \hspace{-.25ex}
 C_{ee}(\cdot) \hspace{-.75ex}\int_s^\cdot \hspace{-1.1ex} U_e(\cdot,r)B_e(r)\wext(r)dr \\
&=  (I-\F_0^s (I-\Delta(\cdot)D_{ee}^\infty)\Delta(\cdot))\inv 
C_{ee}(\cdot)\overline{\mathbb{K}}_s B_e(\cdot)\wext\\
&=  (I-\F_0^s (I-\Delta(\cdot)D_{ee}^\infty)\Delta(\cdot))\inv 
\F_0^s P_{in}\wext.
}
A direct computation shows that $D_{ee}^\infty P_{in} =  D_{ee}^\infty Q_1(\cdot) P_{in} $. Using  $\F_0^s=Q_2(\cdot)[\F_{ee}^\infty - D_{ee}^\infty]Q_1(\cdot)$
with $Q_1(\cdot) = (I-\Delta(\cdot)D_{ee}^\infty)\inv$ and $Q_2(\cdot) = (I-D_{ee}^\infty\Delta(\cdot))\inv$
and denoting $u=P_{in}\wext$ for brevity
we get
\eq{
\F_{ee}^s \wext
&=  (I-Q_2(\cdot)[\F_{ee}^\infty - D_{ee}^\infty]\Delta(\cdot))\inv \F_0^s u + D_{ee}^\infty u\\
&=  (I -\F_{ee}^\infty \Delta(\cdot) )\inv [\F_{ee}^\infty - D_{ee}^\infty]Q_1(\cdot) u+ D_{ee}^\infty Q_1(\cdot) u\\
&=  (I -\F_{ee}^\infty \Delta(\cdot) )\inv \bigl[\F_{ee}^\infty  -\F_{ee}^\infty \Delta(\cdot)D_{ee}^\infty]Q_1(\cdot) u \\
&=  (I -\F_{ee}^\infty \Delta(\cdot) )\inv \F_{ee}^\infty  P_{in}\wext .
}
\potential{
Thus $\F_{ee}^s \wext =  (I -\F_{ee}^\infty \Delta(\cdot) )\inv \F_{ee}^\infty  P_{in}\wext$ on $[s,s+t_0]$ for any $s\geq 0$, $t_0>0$, and $\wext\in \Lploc[2](0,\infty;\C^{n_d+p})$.
}
\end{proof}

\begin{remark}
\label{rem:CLWPremark}
Let $x_0\in X$ and $z_0\in Z$.
If Assumption~\textup{\ref{ass:ContrAss}} is satisfied,
then $(U_e,\Phi_e^{t,s},\Psi_{ee}^s,\F_{ee}^s)_{t\geq s\geq 0}$ is a well-posed nonautonomous system
by
  Theorem~\textup{\ref{thm:CLWellposedness}}.
Thus for every $\wext\in \Lploc[2](0,\infty;\C^{n_d+p})$ the closed-loop state $x_e(t)$ in~\eqref{eq:CLSstateoutput} satisfies $x_e(\cdot)\in \Lploc[2](0,\infty;X_e)$~\textup{\citel{Sch02}{Lem.~3.2}}. Since $K(\cdot)\in\Lp[\infty](0,\infty;\Lin(Z,\C^m))$, for such $\wext(\cdot)$ we also have $u(\cdot)\in\Lploc[2](0,\infty;\C^p)$. 
Theferore the properties of regular linear systems imply that 
if $\wdist(t)$ is as in~\eqref{eq:yrefwdist}, 
 then~\eqref{eq:RLSplant} has a well-defined mild state $x(t)$ satisfying $x(t)\in \Dom(\CL)$ for a.e. $t\geq 0$, and the output $y(t)$ is determined by~\eqref{eq:RLSplantoutput} for a.e. $t\geq 0$.
\end{remark}

\section{Regulation with Converging Controllers}
\label{sec:convfreq}

In this section we introduce general results on output regulation with a nonautonomous controller $(\mc{G}_1(t),\mc{G}_2(t),K(t))$ satisfying Assumption~\ref{ass:ContrAss}.
Our first main result in Theorem~\ref{thm:ORPmain} is applicable in the situation where the controller parameters have well-defined asymptotic limits in the sense that 
\eq{
\begin{cases}
\norm{\Delta_{\mc{G}_1}(t)}\to 0, \\
\norm{\mc{G}_2(t)-\mc{G}_2^\infty}\to 0\\ \norm{K(t)-K^\infty}\to 0
\end{cases}
\qquad \quad \mbox{as} \quad  t\to\infty
}
 for some $\mc{G}_2^\infty\in \Lin(\C^p,Z)$ and $K^\infty\in \Lin(Z,\C^m)$\footnote{Recall 
that 
 ``$\norm{f(t)}\to 0$'' 
for
$f\in \Lp[\infty]$
means $\essup_{s\geq t}\norm{f(s)}\to 0$.
}.
Our second main result in Theorem~\ref{thm:ORPmainNonconv} considers a more general situation where the above norms become small as $t\to\infty$ but do not necessarily converge to zero.
The main condition 
in our results
is that
 the part $\mc{G}_1^\infty$
of $\mc{G}_1(\cdot)$
in Assumption~\ref{ass:ContrAss}
has \emph{an internal model} of the frequencies 
$(\gw_k)_{k=0}^q$ 
(with $\gw_0=0$)
 of $\wdist(t)$ and $\yref(t)$ in~\eqref{eq:yrefwdist}  in the following sense.

\begin{definition}[{\cite[Def. 6.1]{PauPoh10}}]
  \label{def:pcopy}
  The operator $\mc{G}_1^{\infty}$
 \emph{has an internal model of} 
$(\gw_k)_{k=0}^q$ 
if
  $\dim \ker(\pm i\gw_k-\mc{G}_1^\infty)\geq p$ for all $k\in \List[0]{q}$,
where 
 $p\in\N$ is the number of outputs of~\eqref{eq:RLSplant}.
\end{definition}

Our first result states that if the controller parameters converge,
 if $\mc{G}_1^\infty$ has
an internal model of the frequencies 
of $\yref(t)$ and $\wdist(t)$ and if the closed-loop system is exponentially stable, then the controller achieves output regulation.
Exponential stability of $U_e(t,s)$ means that there exists $M,\ga>0$ such that $\norm{U_e(t,s)}\leq Me^{-\ga (t-s)}$ for $t\geq s\geq 0$.

\begin{theorem}
  \label{thm:ORPmain}
  Assume
 $\yref(t)$ and $\wdist(t)$ in~\eqref{eq:yrefwdist} and the initial states $x_0\in X$ and $z_0\in Z$
 are such that
there exist
$\mc{G}_1(\cdot)$, $\mc{G}_2(\cdot)$ and $K(\cdot)$ satisfying Assumption~\textup{\ref{ass:ContrAss}}, and for some $\mc{G}_2^\infty\in \Lin(\C^p,Z)$ and $K^\infty\in \Lin(Z,\C^m)$ 
\eq{
\gd_{\mc{G}}(t):=
\max\set{\norm{\Delta_{\mc{G}_1}(t)},\norm{\mc{G}_2(t)-\mc{G}_2^\infty},\norm{K(t)-K^\infty}}
}
satisfies $\gd_{\mc{G}}(t)\to 0$ as $t\to\infty$.
If 
$U_e(t,s)$ is exponentially stable and $\mc{G}_1^\infty$
if has an internal model of
$(\gw_k)_{k=0}^q$ 
 in~\eqref{eq:yrefwdist},
then
\eqn{
\label{eq:ORPproperty}
\int_t^{t+1} \norm{y(s)-\yref(s)}^2ds\to 0, \qquad \mbox{as} \quad t\to\infty
}
If $\essup_{t\geq 0}e^{\ga t}\gd_{\mc{G}}(t)<\infty$ for some $\ga>0$, 
then there exists 
 $\ga_e>0$ such that 
$t\mapsto e^{\ga_e t}(y(t)-\yref(t))\in \Lp[2](0,\infty;Y)$.
\end{theorem}

In Theorem~\ref{thm:ORPmain} the controller parameters and their limits
 are allowed to depend on the initial states of the system and the controller and of $\wdist(t)$ and $\yref(t)$. 
This possibility is motivated by our controller design for output regulation
 with unknown frequencies in Section~\ref{sec:ContrDesign}.
If $(\mc{G}_1(t),\mc{G}_2(t),K(t))$ are independent of the initial states and $\wdist(t)$ and $\yref(t)$, the claims of Theorem~\ref{thm:ORPmain} (and Theorem~\ref{thm:ORPmainNonconv}) hold for all $x_0\in X$, $z_0\in Z$, $\wdist(t)$ and $\yref(t)$.

The proof of Theorem~\ref{thm:ORPmain} utilises the feedback structure introduced in Section~\ref{sec:CLwellposedness}. To this end 
we
use the notation in Section~\ref{sec:CLwellposedness} and 
in particular denote
 $\Delta_{\mc{G}_2}(\cdot) = \mc{G}_2(\cdot)-\mc{G}_2^\infty\in \Lp[\infty](0,\infty;\Lin(\C^p,Z))$
 and 
 $\Delta_K(\cdot) = K(\cdot)-K^\infty\in \Lp[\infty](0,\infty;\Lin(Z,\C^m))$.
The ``if''-part of the following lemma also follows from~\citel{Sch02}{Thm.~5.6} 
\potential{
(see also~\citel{ClaLat00}{Sec.~4}).%
}%

\begin{lemma}
  \label{lem:EFvsSGstabRLS}
Let $(\mc{G}_1(t),\mc{G}_2(t),K(t))$ satisfy Assumption~\textup{\ref{ass:ContrAss}} and assume $\mc{G}_2^\infty\in \Lin(\C^p,Z)$ and $K^\infty\in \Lin(Z,\C^m)$ are such that $\gd_{\mc{G}}(t)$ in Theorem~\textup{\ref{thm:ORPmain}} satisfies $\gd_{\mc{G}}(t)\to 0$ as $t\to \infty$.
Then 
 $U_e(t,s)$ is exponentially stable if and only if the semigroup $T_e(t)$ generated by $A_e^\infty$ is exponentially stable.
\end{lemma}

\begin{proof}
In the notation of Section~\ref{sec:CLwellposedness}, a direct computation shows 
$ M_D:=\essup_{s\geq 0}\norm{(I-D_{ee}^\infty\Delta(\cdot))\inv}< \infty$.
 Moreover,
Theorem~\ref{thm:CLWellposedness} implies that for every $t_0>0$ there exists  $\gg(t_0)>0$ such that
$\sup_{s\geq 0}\norm{C_{ee}^\infty U_e(s+\cdot,s)x}_{\Lp[2](s,s+t_0)}$ $\leq \gg(t_0)\norm{x}$  for all $x\in X_e$ and $s\geq 0$. Because of this, the integral equation in Theorem~\ref{thm:CLWellposedness} together with the admissibility of $B_{ee}^\infty$ for the semigroup $T_e(t)$ imply that for any fixed $t_0>0$ there exists $M_{t_0}>0$ such that
\eq{
\norm{U_e(s+t_0,s)x-T_e(t_0)x}
&\leq M_{t_0}M_D \norm{\Delta(\cdot)}_{\Lp[\infty](s,\infty)}  \norm{C_{ee}^\infty U_e(\cdot,s)x}_{\Lp[2](s,s+t_0)}\\
&\leq M_{t_0}M_D \norm{\Delta(\cdot)}_{\Lp[\infty](s,\infty)}  \gg(t_0) \norm{x}
}
 for all $s\geq 0$ and $x\in X_e$.
Since $\norm{\Delta(\cdot)}_{\Lp[\infty](s,\infty)}\to 0$ as $s\to\infty$ by assumption, we  have $\norm{U_e(s+t_0,s)-T_e(t_0)}\to 0$ as $s\to\infty$ for every fixed $t_0>0$. 
The claim can now be verified by 
 choosing sufficiently large values of $t_0>0$.
\end{proof}

\begin{proof}[Proof of Theorem~\textup{\ref{thm:ORPmain}}]
Let $x_0\in X$, $z_0\in Z$, $\wdist(t)$ and $\yref(t)$ satisfy the assumptions of the theorem.
According to~\eqref{eq:CLSstateoutput} the regulation error $e(t) = y(t)-\yref(t)$ can be expressed using the output map $\Psi_e^s$ and input-output maps $\F_e^s$ and $\F_e^\infty= P_{out}\F_{ee}^\infty P_{in}$ as 
\eq{
 e(t) &= (\Psi_e^0 x_{e0})(t)+ (\F_e^0 \wext)(t)
\\
&
= (\Psi_e^0 x_{e0})(t)
+ (\F_e^\infty \wext)(t)
+ \bigl[(\F_e^0 \wext)(t)
- (\F_e^\infty \wext)(t)\bigr].
}
If define $B_e^\infty = B_{ee}^\infty P_{in}\in \Lin(\C^{n_d+p},X_e)$ and $C_e^\infty = P_{out}C_{ee}^\infty: \Dom(C_{ee}^\infty)\subset X_e\to \C^p$, then $\F_e^\infty$ is the input-output map of the regular linear system $(A_e^\infty,B_e^\infty,C_e^\infty,D_e)$. A direct computation using~\eqref{eq:CLlimitoperators} shows that this system
has the form of the closed-loop system
 in~\citel{Pau16a}{Sec.~II}
 corresponding to the regular linear system $(A,[B,B_d],C,[D,D_d])$ and the autonomous controller $(\mc{G}_1^\infty,\mc{G}_2^\infty,K^\infty)$. In particular, $(\F_e^\infty \wext)(t)$ is the regulation error corresponding to zero initial states of the system and the controller and the signals $\wdist(t)$ and $\yref(t)$. Since
$\mc{G}_1^\infty$ has an internal model of 
$(\gw_k)_{k=0}^q$
 and 
  $T_e(t)$ generated by $A_e^\infty$ is stable by Lemma~\ref{lem:EFvsSGstabRLS}, we have from~\citel{Pau16a}{Thm.~7} that $t\mapsto e^{\ga_1 t}\norm{(\F_e^\infty \wext)(t)}\in\Lp[2](0,\infty)$ for some $\ga_1>0$.
Since 
 $U_e(t,s)$ is exponentially stable, there exists $\ga_2>0$ such that 
$ t\mapsto e^{\ga_2 t}\norm{(\Psi_e^0 x_{e0})(t)} \in\Lp[2](0,\infty)$.

It remains to analyse the term $\F_e^0 \wext -\F_e^\infty \wext$.
Note that
  $\sup_{t\geq 0}\norm{\wext}_{\Lp[2](t,t+1)}$ $<\infty$.
For $\Delta(\cdot)$ in~\eqref{eq:Deltat} 
we have $\essup_{t\leq s\leq t+1}\norm{\Delta(s)}\to 0$ as $t\to\infty$ if
 $\gd_{\mc{G}}(t)\to 0$ as $t\to \infty$.
Theorem~\ref{thm:CLWellposedness} implies
$ \F_{ee}^s = (I-\F_{ee}^\infty\Delta(\cdot))\inv  \F_{ee}^\infty P_{in}$ and%
\begin{subequations}
\label{eq:IOmapdifference}
\eqn{
\label{eq:IOmapdifference1}
\F_e^0 -\F_e^\infty
 &= P_{out}\bigl[(I-\F_{ee}^\infty\Delta(\cdot))\inv - I\bigr] \F_{ee}^\infty P_{in}\\
\label{eq:IOmapdifference2}
 &= P_{out}\F_{ee}^\infty \Delta(\cdot)(I-\F_{ee}^\infty\Delta(\cdot))\inv  \F_{ee}^\infty P_{in}\\
 &= P_{out}\F_{ee}^\infty \Delta(\cdot)  \F_{ee}^0.
}
\end{subequations}
As shown in the proof of Theorem~\ref{thm:CLWellposedness}, $\F_{ee}^s$ 
 is an input-output map of a nonautonomous well-posed system
 with an exponentially stable evolution family $U_e(t,s)$. 
Thus Lemma~\ref{lem:IOmapconvproperties}(a) implies $\sup_{t\geq 0}\norm{\F_{ee}^0 \wext}_{\Lp[2](t,t+1)}<\infty$ and
\eq{
\MoveEqLeft[1]\norm{\Delta(\cdot) \F_{ee}^0  \wext}_{\Lp[2](t,t+1)}
\leq \norm{\Delta(\cdot)}_{\Lp[\infty](t,t+1)}\norm{\F_{ee}^0 \wext}_{\Lp[2](t,t+1)}
\to 0
}
as $t\to\infty$.
Moreover, Lemma~\ref{lem:EFvsSGstabRLS} implies 
that the regular linear system $(A_e^\infty,B_{ee}^\infty,C_{ee}^\infty,D_{ee}^\infty)$ with extended input-output map $\F_{ee}^\infty$ is exponentially stable.
Thus Lemma~\ref{lem:IOmapconvproperties}(c)
applied to this autonomous system and
$u=\Delta(\cdot) \F_{ee}^0\wext$
together with~\eqref{eq:IOmapdifference} show that
\ieq{
\norm{\F_e^0 \wext -\F_e^\infty \wext}_{\Lp[2](t,t+1)}
\to 0
}
as $t\to\infty$.
This completes the proof of~\eqref{eq:ORPproperty}.

Finally, let $\ga>0$ be such that
 $\essup_{t\geq 0}e^{\ga t}\gd_{\mc{G}}(t)<\infty$.
Then we have $\sup_{t\geq 0} e^{\ga t} \norm{\Delta(\cdot)}_{\Lp[\infty](t,t+1)}<\infty$
and
since
$\sup_{t\geq 0}\norm{\F_{ee}^0 \wext}_{\Lp[2](t,t+1)}<\infty$, we can estimate
\eq{
\sup_{t\geq 0} e^{\ga t} \norm{\Delta(\cdot)  \F_{ee}^0  \wext}_{\Lp[2](t,t+1)}
&\leq\sup_{t\geq 0} e^{\ga t} \norm{\Delta(\cdot)}_{\Lp[\infty](t,t+1)}\norm{ \F_{ee}^0  \wext}_{\Lp[2](t,t+1)}
<\infty.
}
Lemma~\ref{lem:IOmapconvproperties}(d) 
for $(A_e^\infty,B_{ee}^\infty,C_{ee}^\infty,D_{ee}^\infty)$ and
 $u=\Delta(\cdot)  \F_{ee}^0 \wext$
together with~\eqref{eq:IOmapdifference} imply
$t\mapsto e^{\ga_0 t}\norm{(\F_e^0 \wext)(t)-(\F_e^\infty \wext)(t)}\in\Lp[2](0,\infty)$ for some
  $\ga_0>0$. This completes the proof.
\end{proof}

The following result generalises Theorem~\ref{thm:ORPmain} to the situation where the parameters of the controller do not necessarily converge as $t\to \infty$, 
or where the limit of $\mc{G}_1(t)$ has an internal model of frequencies which are only close to
$(\gw_k)_{k=0}^q$.
\potential{%
In these cases the asymptotic tracking error will be small provided that the asymptotic error in the frequencies is sufficiently small.
}

\begin{theorem}
  \label{thm:ORPmainNonconv}
  Assume that
 $x_0\in X$, $z_0\in Z$,
$\yref(t)$, and $\wdist(t)$ in~\eqref{eq:yrefwdist} 
 are such that
there exist
$\mc{G}_1(\cdot)$, $\mc{G}_2(\cdot)$ and $K(\cdot)$ satisfying Assumption~\textup{\ref{ass:ContrAss}} and 
$U_e(t,s)$ is exponentially stable. 
Moreover, assume 
 $\mc{G}_1^\infty$ has an internal model of
$(\gw_k)_{k=0}^q$
 and 
$\mc{G}_2^\infty\in \Lin(\C^p,Z)$ and $K^\infty\in \Lin(Z,\C^m)$ are such that $T_e(t)$ is exponentially stable.
Define
\eq{
\gd_{\mc{G}}(t):=
\max\set{\norm{\Delta_{\mc{G}_1}(t)},\norm{\mc{G}_2(t)-\mc{G}_2^\infty},\norm{K(t)-K^\infty}}
}
There exist $\Merr,\gd_0>0$
 depending only on $(A,B,B_d,C,D,D_d)$ and $(\mc{G}_1^\infty,\mc{G}_{11}^\infty,$ $\mc{G}_2^\infty,K^\infty)$
such that
\eq{
\MoveEqLeft[3]\limsup_{t\to\infty}\int_t^{t+1} \norm{y(s)-\yref(s)}^2ds
\\
&
\leq \Merr 
\norm{[\wdist(\cdot),\yref(\cdot)]^T}_\infty^2
\limsup_{t\to\infty}\norm{\gd_{\G}(\cdot)}_{\Lp[\infty](t,\infty)}^2
}
provided that
$\limsup_{t\to\infty}\norm{\gd_{\G}(\cdot)}_{\Lp[\infty](t,\infty)}< \gd_0$.
\end{theorem}

\begin{proof}
As shown in the proof of Theorem~\ref{thm:ORPmain}, we have
\ieq{
 e(\cdot) = e_0(\cdot)
+ (\F_e^0 \wext
- \F_e^\infty \wext),
}
where  $w_e(t)=[\wdist(t),\yref(t)]^T$ and
$ t\mapsto e^{\ga_0 t}\norm{e_0(t)} \in\Lp[2](0,\infty)$ for some $\ga_0>0$.
Moreover, the 
identity~\eqref{eq:IOmapdifference}
 and Lemma~\ref{lem:IOmapconvproperties}(b) imply that there exists $M_1>0$ depending only on $(A_e^\infty,B_{ee}^\infty,C_{ee}^\infty,D_{ee}^\infty)$ such that 
\eq{
\limsup_{t\to\infty}\norm{e(\cdot)}_{\Lp[2](t,t+1)}
&\leq \limsup_{t\to\infty}\;\norm{\F_e^0 \wext - \F_e^\infty \wext}_{\Lp[2](t,t+1)}\\
&\leq M_1\limsup_{t\to\infty}\Bigl(\norm{\Delta(\cdot)}_{\Lp[\infty](t,t+1)}\norm{\F_{ee}^0  \wext}_{\Lp[2](t,t+1)}\Bigr),
}
where
$\F_{ee}^0= (I-\F_{ee}^\infty\Delta(\cdot))\inv  \F_{ee}^\infty P_{in}$.
Since $\F_{ee}^\infty$ is the extended input-output map of the regular linear system $(A_e^\infty,B_{ee}^\infty,C_{ee}^\infty,D_{ee}^\infty)$, by Lemma~\ref{lem:IOmapconvproperties}(a) there exists $M_0>0$ such that $\sup_{\tau\geq 0}\norm{\F_{ee}^\infty u}_{\Lp[2](\tau,\tau+1)}\leq M_0 \sup_{\tau\geq 0}\norm{ u}_{\Lp[2](\tau,\tau+1)}$ for all $u\in \Lploc[2](0,\infty;U_{ee})$.
We define
\eq{
\gd_0 = 
2^{-3/2} \min \set{\norm{\F_{ee}^\infty}\inv, M_0\inv}
}
and assume $\essup_{t\geq t_0}\gd_{\mc{G}}(t)\leq \gd_0$ for some fixed $t_0>0$.
Since 
the definition of $\Delta(t)$ in~\eqref{eq:Deltat} implies that $\norm{\Delta(t)}\leq \sqrt{2}\gd_{\mc{G}}(t)$ for a.e. $t\geq 0$, we have 
$\norm{\F_{ee}^\infty}\norm{\Delta}_{\Lp[\infty](t_0,\infty)}\leq 1/2$ and 
$M_0\norm{\Delta}_{\Lp[\infty](t_0,\infty)}\leq 1/2$.

By Theorem~\ref{thm:CLWellposedness},
$(U_e,\Phi_e^{t,s},\Psi_{ee}^s,\F_{ee}^s)_{t\geq s\geq 0}$ is an exponentially stable nonautonomous well-posed system.
We have from~\citel{Sch02}{Def.~3.6} that for a.e. $t\geq t_0$
\eqn{
\label{eq:ORPnonconvFeesplit}
(\F_{ee}^0 \wext)(t) = (\Psi_{ee}^{t_0}\Phi_e^{t_0,0}\wext)(t) + (\F_{ee}^0w_{t_0})(t),
}
where $w_{t_0}: [0,\infty)\to\C^p$ is defined so that $w_{t_0}(t)=0$ for $t\in[0,t_0)$ and $w_{t_0}(t)=\wext(t)$ for $t\geq t_0$.
Define $\Delta_{t_0}\in \Lp[\infty](0,\infty;\Lin(Y_{ee},U_{ee}))$ so that $\Delta_{t_0}(t)=0$ for $t\in [0,t_0)$ and $\Delta_{t_0}(t)=\Delta(t)$ for a.e. $t\geq t_0$. Then
$\norm{\Delta_{t_0}}_{\Lp[\infty]}=\norm{\Delta}_{\Lp[\infty](t_0,\infty)}$, and the properties  
$(\F_{ee}^\infty P_{in}w_{t_0})(t)=0$ for $t\in [0,t_0]$ and $\norm{\F_{ee}^\infty}\norm{\Delta_{t_0}(\cdot)}_{\Lp[\infty]}\leq 1/2$ imply
\eqn{
\label{eq:ORPNCFee0formula}
\F_{ee}^0 w_{t_0} 
&= (I-\F_{ee}^\infty\Delta_{t_0}(\cdot))\inv  \F_{ee}^\infty P_{in} w_{t_0}
= \sum_{n=0}^\infty (\F_{ee}^\infty\Delta_{t_0}(\cdot))^n  \F_{ee}^\infty P_{in} w_{t_0}
}%
with convergence
 in $\Lploc[2](0,\infty;Y_{ee})$.
The choice of $M_0>0$ 
and 
 $\norm{\wext}_{\Lp[2](t,t+1)} $ $\leq \norm{\wext}_\infty$ for $t\geq 0$ 
imply that for all $n\in \N$
\eq{
\MoveEqLeft\sup_{t\geq 0}\; \norm{(\F_{ee}^\infty\Delta_{t_0}(\cdot))^n  \F_{ee}^\infty P_{in} w_{t_0}}_{\Lp[2](t,t+1)}
\\
&\leq (M_0\norm{\Delta_{t_0}}_{\Lp[\infty]})^n \sup_{t\geq 0}\;\norm{\F_{ee}^\infty P_{in} w_{t_0}}_{\Lp[2](t,t+1)}\\
&\leq 2^{-n} M_0  \sup_{t\geq 0}\;\norm{  P_{in}w_{t_0}}_{\Lp[2](t,t+1)}
\leq 2^{-n} M_0  \norm{\wext}_\infty.
}
Thus~\eqref{eq:ORPNCFee0formula} implies $\sup_{t\geq 0}\norm{\F_{ee}^0 w_{t_0}}_{\Lp[2](t,t+1)}\leq 2M_0
\norm{\wext}_\infty
 $.
Since
we have
 $\norm{\Psi_{ee}^{t_0}\Phi_e^{t_0,0}\wext}_{\Lp[2](t,t+1)}\to 0$ as $t\to\infty$ and
since
 $\norm{\Delta(t)}\leq \sqrt{2}\gd_{\G}(t)$
for a.e. $t\geq 0$,
 equation~\eqref{eq:ORPnonconvFeesplit} finally implies%
\eq{
\limsup_{t\to\infty}\;\norm{e(\cdot)}_{\Lp[2](t,t+1)}
&\leq  M_1\limsup_{t\to\infty}\Bigl(\norm{\Delta(\cdot)}_{\Lp[\infty](t,t+1)}
\norm{\F_{ee}^0  w_{t_0}}_{\Lp[2](t,t+1)}
\Bigr)\\
&
\potential{
\leq  2M_0M_1 \norm{w_e}_\infty
\limsup_{t\to\infty}\;\norm{\Delta(\cdot)}_{\Lp[\infty](t,t+1)}
}
\\
&\leq  2 \sqrt{2} M_0M_1 \norm{w_e}_\infty
\lim_{t\to\infty}\;\norm{\gd_{\G}(\cdot)}_{\Lp[\infty](t,\infty)}.
\qedhere
}
\end{proof}

\begin{remark}
\label{rem:ORPnonconvDeltaDependence}
The proof of 
 Theorem~\textup{\ref{thm:ORPmainNonconv}} 
shows that
 $\Merr,\gd_0>0$
depend
on
the norm $\norm{\F_{ee}^\infty}$ of the extended input-output map
of 
the autonomous system $(A_e^\infty,B_{ee}^\infty,C_{ee}^\infty,D_{ee}^\infty)$
and
 on 
 $M_0,M_1>0$  
in Lemma~\textup{\ref{lem:IOmapconvproperties}(a)--(b)} corresponding to this system.
 Moreover, 
 the proof of Lemma~\textup{\ref{lem:IOmapconvproperties}},
implies that
$M_0,M_1>0$ are determined by
 constants $M_e,\ga_e>0$ such that $\norm{T_e(t)}\leq M_ee^{-\ga_et}$ for $t\geq 0$ and by 
 upper bounds for the norms of the input, output, and input--output map of
 $(A_e^\infty,B_{ee}^\infty,C_{ee}^\infty,D_{ee}^\infty)$ 
 on 
the time interval
 $[0,1]$.
\end{remark}

\begin{remark}
\label{rem:Gconds}
By~\textup{\citel{Pau16a}{Thm.~7}}
the internal model of $\mc{G}_1^\infty$
in Theorems~\textup{\ref{thm:ORPmain}} and~\textup{\ref{thm:ORPmainNonconv}}
 can
 be replaced with the conditions 
\begin{subequations}
\label{eq:Gconds}
\eqn{
\hspace{-1.5ex}\ran(\pm i\gw_k-\mc{G}_1^\infty)\cap \ran(\mc{G}_2^\infty)&=\set{0} \quad \forall k\in \List[0]{q}\\
\ker(\mc{G}_2^\infty)&=\set{0}.
}
\end{subequations}
\end{remark}

\section{Controller Design for
 Output Regulation with
 Unknown Frequencies}
\label{sec:ContrDesign}

In this section we introduce a controller for output regulation of $\yref(t)$ and $\wdist(t)$ with unknown frequencies, amplitudes, and phases.
Our controller contains a time-varying internal model of \textit{estimates} $(\hat{\gw}_k(t))_{k=1}^q$ of  $(\gw_k)_{k=1}^q$ in~\eqref{eq:yrefwdist} and an observer-based part for achieving closed-loop stability.
The frequency estimates are formed based on an auxiliary output $\yaux(t)$ of the controller which contains the information on $(\gw_k)_k$ but 
is by design independent of the time-varying parts of the controller.
Therefore our controller
 can be combined with any estimation method which can asymptotically estimate the frequencies $(\gw_k)_k$ from $\yaux(t)$. 

We introduce the general controller structure  in Section~\ref{sec:ObsController} and 
present general conditions for output regulation
in the situation where $(\hat{\gw}_k(t))_{k=1}^q$ converge to  $(\gw_k)_k$ in~\eqref{eq:yrefwdist}.
 We analyse the structure of $\yaux(t)$ in Lemma~\ref{lem:yauxindep} and  
in Remark~\ref{rem:FreqEstReview} 
we list selected methods which can be used to
estimate the frequencies based on $\yaux(t)$.
In Section~\ref{sec:TuningAlgorithm} we present a tuning algorithm for constructing the
 controller parameters in order to achieve closed-loop stability and output regulation.
Finally, in
Section~\ref{sec:ContrRobAnalysis} we analyse the robustness properties of our controller.
We make the following standing assumptions.

\begin{assumption}
\label{ass:RLSstabilizability}
There exist $\tilde{K}\in \Lin(X,\C^m)$ and $L\in \Lin(\C^p,X)$ such that the semigroups generated by $A+B\tilde{K}: \Dom(A+B\tilde{K})\subset X\to X$ with domain $\Dom(A+B\tilde{K})=\setm{x\in X}{Ax+B\tilde{K}x\in X}$ and $A+LC: \Dom(A)\subset X\to X$ are exponentially stable.
\end{assumption}

\begin{definition}[{\citel{NatGil14}{Def.~V.1}}]
\label{def:TransmissionZero}
The point 
 $i\gl_0\in i\R$
 \emph{is a transmission zero of} $(A,B,C,D)$ if
$P_{\tilde{K}}(i\gl_0)$ is not surjective, 
where
$\tilde{K}\in \Lin(X,\C^m)$ is such that $i\gl_0\in\rho(A+B \tilde{K})$ and $P_{\tilde{K}}(\gl)$ is the transfer function of
$(A+B\tilde{K},B,C+D\tilde{K},D)$.
\end{definition}

\begin{assumption}
  \label{ass:TZass}
Assume $\yref(t)$ and $\wdist(t)$ are of the form~\eqref{eq:yrefwdist} with $0=\gw_0<\gw_1<\ldots<\gw_q$ and
 $(A,B,C,D)$ does not have transmission zeros at $\set{0}\cup\set{\pm i\gw_k}_{k=1}^q$.
\end{assumption}

\subsection{Controller with a time-varying internal model}
\label{sec:ObsController}

Our error feedback controller has the form
\begin{subequations}
\label{eq:ContrFull}
\eqn{
  \dot{z}(t)&= \mc{G}_1(t)z(t)+\mc{G}_2e(t), \quad z(0)=z_0\in Z\\
  u(t)&= K(t)z(t)\\
\label{eq:ContrFullyaux}
\yaux(t)&= \Kaux(t)  z(t) + e(t)
}
\end{subequations}
where
$e(t) = y(t)-\yref(t)$ is the regulation error.
The controller structure in  Definition~\ref{def:ObsContr} generalises the autonomous robust controller in~\citel{HamPoh10}{Sec.~7} and the adaptive internal model based controller scheme in~\citel{CarGal16}{Sec.~4}, where a separate ``residual generator'' system was used to construct $\yaux(t)$.

\begin{definition}
\label{def:ObsContr}
The controller $(\mc{G}_1(t),\mc{G}_2,K(t))$ 
on 
 $Z=Z_0\times X$ with $Z_0=\C^{p(2q +1)}$ is defined
by choosing $L\in \Lin(\C^p,X)$  so that $A+LC$ generates an exponentially stable semigroup $T_L(t)$, 
$K(\cdot)=[K_1(\cdot),K_2(\cdot)]\in \Lp[\infty](0,\infty;\Lin(Z,\C^m))$,  and
\eq{
\MoveEqLeft[3.5]
\mathcal{G}_1(t)=
\pmat{ G_1(t)&0\\
 (B+LD)K_1(t)&A+L\CL+(B+LD)K_2(t) }
\\
\Dom(\mc{G}_1(t))&=\setm{\pmatsmall{z_0\\x}\in Z_0\times \Dom(\CL)}{Ax+BK(t)\pmatsmall{z_0\\x}\in X}
}
\eq{
\mathcal{G}_2&=
\pmat{
G_2\\
-L
},
\; \Kaux(t) = [-DK_1(t),-\CL-DK_2(t)]
}
 with $\Dom(\Kaux(t))=Z_0\times \Dom(\CL)$.
Finally, we define
\eq{
  G_1(t)
&= \diag(0_p, \hat{\gw}_1(t)\Omega_p,\ldots, \hat{\gw}_{q }(t)\Omega_p)\in \Lin(Z_0),\\
  G_2 &= \bigl[I_p,I_p,0_p,I_p,0_p,\ldots,I_p,0_p\bigr]^T \in \R^{p(2q +1)\times p}
}
with $\Omega_p = \pmatsmall{0_{p}& I_p\\-I_p&0_p}$,
where $0_p,I_p\in \R^{p\times p}$ are the 
 zero and identity matrices and
$\hat{\gw}_k(\cdot)\in \Lp[\infty](0,\infty;\R_+)$ for all $k$.
\end{definition}

The function $G_1(\cdot)$ is the \emph{time-varying internal model} which contains the 
estimates  $(\hat{\gw}_k(t))_{k=1}^q$ of the frequencies in $\wdist(t)$ and $\yref(t)$.
By construction, for every $t\geq 0$ the pair $(G_1(t),G_2)$ is controllable if the values $(\hat{\gw}_k(t))_{k=1}^q$ are distinct and nonzero.
If 
 $\abs{\hat{\gw}_k(t)- \gw_k}\to0 $ as $t\to\infty$ for all $k$, then $\norm{G_1(t)- G_1^\infty}\to 0$ as $t\to\infty$ where $G_1^\infty \in \Lin(Z_0) $ is defined by replacing $(\hat{\gw}_k(t))_k$ in $G_1(t)$ with $(\gw_k)_k$.

For any $G_1^\infty\in \Lin(Z_0)$ and 
$K^\infty=[K_1^\infty,K_2^\infty]\in \Lin(Z,\C^m)$ 
we can define $\Delta_G(t)=G_1(t)-G_1^\infty$, $\Delta_K(t)=K(t)-K^\infty$ and
\begin{subequations}
\label{eq:ObsContrLimit}
\eqn{
\hspace{-.8ex}\mc{G}_1^\infty&=\pmat{G_1^\infty&0\\(B+LD) K_1^\infty&A+L\CL+(B+LD)K_2^\infty}\\
\hspace{-.8ex}\mc{G}_{11}^\infty&=\pmat{I&0\\0&B+LD}, \qquad 
\Delta_{\mc{G}_1}(t)=\pmat{[\Delta_G(t),0]\\\Delta_K(t)}.
}
\end{subequations}
The feedback theory for regular linear systems in~\cite{Wei94} implies that 
 $\mc{G}_1^\infty$
with domain $\Dom(\mc{G}_1^\infty)=
\setm{\pmatsmall{z_0\\x}\in Z_0\times \Dom(\CL)}{Ax+BK^\infty\pmatsmall{z_0\\x}\in X}
$ generates a strongly continuous semigroup on $Z$ and that $\mc{G}_{11}^\infty$ is an admissible input operator for this semigroup. Since $\mc{G}_1(t)=\mc{G}_1^\infty + \mc{G}_{11}^\infty \Delta_{\mc{G}_1}(t)$ and $\Delta_{\mc{G}_1}(\cdot)\in \Lp[\infty](0,\infty;\Lin(Z,Z_0\times U))$, the controller in Definition~\ref{def:ObsContr} satisfies Assumption~\ref{ass:ContrAss}.
Therefore the well-defined mild closed-loop state $x_e(t)$ and regulation error $e(t)$ are guaranteed by
Theorem~\ref{thm:CLWellposedness}.

Our first result shows that if the frequency estimates $(\hat{\gw}_k(t))_k$ converge to the true frequencies and if $K(\cdot)$ is such that the semigroup generated by the block operator~\eqref{eq:SGcontrdesign} is stable, then
 the controller achieves output regulation.

\begin{theorem}
  \label{thm:ContrMain}
  Choose $(\mc{G}_1(t),\mc{G}_2,K(t))$ as in Definition~\textup{\ref{def:ObsContr}}.
Assume
that
$\wdist(t)$ and $\yref(t)$ and the intial conditions $x_0\in X$ and $z_0\in Z$
 are such that 
\eq{
\abs{\hat{\gw}_k(t)-\gw_k}\to 0 \quad \mbox{and} \quad
\norm{K(t)-K^\infty}\to 0
}
as $t\to\infty$ 
for all $k$ and for some $K^\infty\in \Lin(Z,\C^m)$.
If 
the semigroup generated by
  \eqn{
    \label{eq:SGcontrdesign}
A_s^\infty+B_sK^\infty:= 
   \pmat{G_1^\infty &G_2\CL\\0&A}+\pmat{G_2D\\B}K^\infty
  }
with domain $\setm{[z_0,x]^T\in Z_0\times \Dom(\CL)}{Ax+BK^\infty \pmatsmall{z_0\\x}\in X}$
is exponentially stable, then
\eq{
\int_t^{t+1} \norm{y(s)-\yref(s)}^2ds\to 0, \qquad \mbox{as} \quad t\to\infty
}
and $U_e(t,s)$ is exponentially stable.
  If  $\essup_{t\geq 0}e^{\ga t}\abs{\hat{\gw}_k(t)- \gw_k}<\infty$ and
$\essup_{t\geq 0}e^{\ga t}\norm{K(t)-K^\infty}<\infty$ for some $\ga>0$ and for all 
$k$, 
then there exists 
 $\ga_e>0$ such that 
$t\mapsto e^{\ga_e t}(y(t)-\yref(t))\in \Lp[2](0,\infty;Y)$.
\end{theorem}

\begin{proof}
Assume that $\wdist(t)$, $\yref(t)$, $x_0\in X$, and $z_0\in Z$ are such the assumptions hold. If we define 
$G_1^\infty \in \Lin(Z_0) $ by replacing $(\hat{\gw}_k(t))_k$ in $G_1(t)$ with $(\gw_k)_k$ and let $\Delta_G(t)=G_1(t)-G_1^\infty$ and $\Delta_K(t)=K(t)-K^\infty$, then $\Delta_{\mc{G}_1}(t)$ in~\eqref{eq:ObsContrLimit} satisfies $\norm{\Delta_{\mc{G}_1}(t)}\to 0$ as $t\to\infty$.
 As shown in the proof of~\citel{Pau16a}{Thm.~15}, the pair $(\mc{G}_1^\infty,\mc{G}_2)$ satisfies the ``$\mc{G}$-conditions''~\eqref{eq:Gconds}.
In view of Remark~\ref{rem:Gconds}, the claims  follow from Theorem~\ref{thm:ORPmain} and Lemma~\ref{lem:EFvsSGstabRLS} once we show that the semigroup $T_e(t)$ generated by $A_e^\infty$ is exponentially stable.

The operator $A_e^\infty$ is exactly the operator $A_e(t)$ with $\mc{G}_1(t)$ and $K(t)$ replaced with $\mc{G}_1^\infty$ and $K^\infty$, respectively, i.e.,
\eq{
  A_e^\infty
  = \pmat{A&BK_1^\infty &BK_2^\infty \\G_2\CL&G_1^\infty +G_2DK_1^\infty &G_2DK_2^\infty \\-L\CL&B K_1^\infty &A+L\CL+BK_2^\infty }.
}
If we define
 $Q_e\in \Lin(X\times Z_0\times X,Z_0\times X\times X)$ by 
\eqn{
\label{eq:Qetransform}
  Q_e =  \pmat{0&I&0\\I&0&0\\-I&0&I}, 
  \quad 
  Q_e\inv =  \pmat{0&I&0\\I&0&0\\0&I&I}, 
}
a direct computation shows that 
\eqn{
\label{eq:QeAeform}
  Q_eA_e^\infty Q_e\inv 
  &=\pmat{A_s^\infty+B_sK^\infty&B_sK_2^\infty\\
  0&A+LC}
}
with 
$\Dom(Q_eA_e^\infty Q_e\inv)= \setm{[x_s,\tilde{x}]^T\in (Z_0\times \Dom(\CL))\times \Dom(A)}{ A_sx_s+B_s(K^\infty x_s+K_2^\infty \tilde{x})\in Z_0\times X}$.
By assumption the semigroups generated by $A+LC: \Dom(A)\subset X\to X$ and $A_s+B_sK^\infty$ are exponentially stable. Moreover, $B_s$ is an admissible input operator for the semigroup generated by $A_s+B_sK^\infty$ by the results in~\citel{Wei94}{Sec.~7}.
Thus  
 the semigroup generated by $Q_eA_e^\infty Q_e\inv$ is exponentially stable and similarity implies the same for $T_e(t)$.
The claims now follow from Theorem~\ref{thm:ORPmain} with Lemma~\ref{lem:EFvsSGstabRLS} and Remark~\ref{rem:Gconds}.
\end{proof}

We conclude this section by
 analysing
 $\yaux(t)$.
Lemma~\ref{lem:yauxindep} in particular shows that $\yaux(t)$ is independent of the time-varying parameters 
$(\hat{\gw}_k(t))_{k=1}^q$ and $K(t)$.
The form of $\yaux(t)$ involves $B_{dL}=[B_d+LD_d,-L]$ and 
the transfer function 
\eqn{
\label{eq:Ptot}
\hspace{-1ex}\Ptot(\gl) = \CL R(\gl,A+L\CL)B_{dL}+ [D_d,-I]
}
of the regular linear system $(A+LC,B_{dL},C,[D_d,-I])$.
In our controller in Section~\ref{sec:TuningAlgorithm} the parameters $G_1(t)$ and $K(t)$ are designed based on $\yaux(t)$, but the properties $G_1(\cdot)\in\Lp[\infty](0,\infty;\Lin(Z_0))$ and $K(\cdot)\in \Lp[\infty](0,\infty;\Lin(Z,\C^m))$  
 are guaranteed by construction, and thus the structure~\eqref{eq:yauxform} of $\yaux(t)$ is determined by Lemma~\ref{lem:yauxindep}.

\begin{lemma}
\label{lem:yauxindep}
Let $x_0\in X$ and $z_0=(z_{10},z_{20})\in Z$ and let $\yref(t)$ and $\wdist(t)$ be as in~\eqref{eq:yrefwdist}.
  Consider the controller $(\mc{G}_1(t),\mc{G}_2,K(t))$ in Definition~\textup{\ref{def:ObsContr}}, let $T_L(t)$ be the semigroup generated by $A+LC$
and
denote $\gw_{-k}:=-\gw_k$ for $k\in \List{q}$. Then 
\eqn{
\label{eq:yauxform}
\yaux(t)=y_0(t)+\sum_{k=-q}^q
e^{i\gw_k t} \Ptot(i\gw_k)  \cextk 
}
for a.e. $t\geq 0$
with $\cextk = \frac{1}{2} [\wdistamplk e^{i\varphi_k},\yrefamplk e^{i\theta_k}]^T$, $\cextk[-k]=\conj{\cextk}$, and%
\eq{
 y_0(t) &= \CL T_L(t)\bigl(x_0-z_{20}- \hspace{-.5ex}\sum_{k=-q}^q \hspace{-.5ex}R(i\gw_k,A+LC)B_{dL}c_e^k\bigr).
}
We have $t\mapsto e^{\ga t}\ytrans(t)\in \Lp[2](0,\infty;\C^p)$  for some $\ga>0$.
Moreover, if 
$C\in \Lin(X,\C^p)$, if
$x_0-z_{20}\in \Dom(A)$, or if $A$ generates an analytic semigroup,
then $y_0(\cdot)$ is continuous and
 $ e^{\ga t}\norm{\ytrans(t)}\to 0$ as $t\to \infty$  for some $\ga>0$.
\end{lemma}

\begin{proof}
Denote
$A_L=A+LC$,
 $B_{dL}=[B_d+LD_d,-L]$, and
\eq{
A_s(t)=\pmat{G_1(t)&G_2\CL\\0&A}.
}
Since Assumption~\ref{ass:ContrAss} holds
the
closed-loop state $x_e(t)=[x(t),z_1(t),z_2(t)]^T$ and the regulation error $e(t)$  in~\eqref{eq:CLSstateoutput} are well-defined by Theorem~\ref{thm:CLWellposedness}.
With $B_s\in \Lin(U,Z_0\times X_{-1})$ in~\eqref{eq:SGcontrdesign} and
$Q_e$ in~\eqref{eq:Qetransform} we have 
\eq{
Q_eB_e(t)
 &\equiv Q_e \pmat{B_d&0\\ G_2 D_d&-G_2\\-LD_d&L} 
= \pmat{ G_2 D_d&-G_2\\B_d&0\\-B_d-LD_d&L}\\
Q_eA_e(t)Q_e\inv 
 & =\pmat{A_s(t)+B_sK(t)&B_sK_2(t)\\
  0&A_L}
}
with
 $\Dom(Q_eA_e(t)Q_e\inv ) =\setm{[x_s,\tilde{x}]^T\in (Z_0\times \Dom(\CL))\times \Dom(A)}{A_s(t)x_s+B_s(K(t)x_s+K_2(t)\tilde{x})\in Z_0\times X}$ 
for a.e. $t\geq 0$.
Consequently $Q_eU_e(t,s)Q_e\inv$ has a block triangular form for all $t\geq s\geq 0$.
Applying the similarity transformation $Q_e$ in~\eqref{eq:Qetransform} to~\eqref{eq:CLSstateoutputstate}  and~\eqref{eq:CLMapsDefPhi} therefore shows that
$\tilde{x}(t) = x(t)-z_2(t)$ is the mild solution of  
\eqn{
\label{eq:yauxformACP}
\dot{\tilde{x}}(t)
= A_L \tilde{x}(t)
 +B_{dL} \wext(t), \qquad \tilde{x}(0)=x_0-z_{20}
}
with $\wext(t)=[\wdist(t),\yref(t)]^T$.
Since 
$B_{dL}\in \Lin(\C^{n_d+p},X)$ and since $C$ is an admissible 
with respect to
  $T_L(t)$ generated by $A_L$, we have $\tilde{x}(t)\in \Dom(\CL)$ for a.e. $t\geq 0$. 
Moreover, $x(t)\in \Dom(\CL)$ for a.e. $t\geq0$ by 
 Theorem~\ref{thm:CLWellposedness}.
Thus $z_2(t)= x(t)-\tilde{x}(t)\in \Dom(\CL)$ for a.e. $t\geq 0$ and 
the formula~\eqref{eq:ContrFullyaux} for $\yaux(t) $ is well-defined for a.e. $t\geq 0$.
Since $\Kaux(t)z(t)=-\CL z_2(t)-DK(t)z(t)$, Remark~\ref{rem:CLWPremark} and~\eqref{eq:RLSplantoutput} imply that 
\eq{
\yaux(t)
&=-\CL z_2(t) - Du(t)+y(t)-\yref(t)
= \CL\tilde{x}(t)+[D_d,\,-I]\wext(t)
}
for a.e. $t\geq 0$.
Thus $\yaux(t)$ is the output of the regular linear system $(A_L,B_{dL},C,[D_d,-I])$ with initial state $\tilde{x}(0)=x_0-z_{20}\in X$ and input $\wext(t)$. When $\tilde{x}(0)=0$ and $\wext(t)=e^{i\gw_k t} w_0$ for some $w_0\in \C^{n_d+p}$, \citel{Sta05book}{Cor.~4.6.13} implies
\eq{
\yaux(t)&=
  e^{i\gw_k t}\Ptot(i\gw_k)w_0
-\CL T_L(t)R(i\gw_k,A_L)B_{dL} w_0 .
}
Finally, linearity implies that for $\tilde{x}(0)=x_0-z_{20}$ and $\wext(t)=[\wdist(t),\yref(t)]^T$ the output $\yaux(t)$ has the form in~\eqref{eq:yauxform} 
with the given $\set{\cextk}_{k=-q}^q$ and $y_0(t)$.
Since $C$ is admissible with respect to the exponentially stable semigroup $T_L(t)$, we have $t\mapsto e^{\ga t}y_0(t)\in \Lp[2](0,\infty;\C^p)$ for some $\ga>0$.

In the last claim, if
 $C\in \Lin(X,\C^p)$, then pointwise convergence of $\ytrans(t)$ follows directly from stability of $T_L(t)$. In the other cases
the property
$x_1:=\sum_{k=-q}^q R(i\gw_k,A_L)B_{dL}c_e^k\in \Dom(A_L)$ and
  $CA_L\inv\in \Lin(X,\C^p)$ 
imply that
\ieq{
\ytrans(t)= CA_L\inv A_L T_L(t)(x_0-z_{20}-x_1)\to 0
}
at an exponential rate as $t\to\infty$.
\end{proof}

\begin{remark}[Methods for Frequency Estimation]
\label{rem:FreqEstReview}
 Multi-frequency estimators 
based on dynamical adaptive observers have been developed by several authors in, 
e.g.,~\cite{ObrCas02,MarTom02,Xia02,CarAst08,AraBob16,WanGri17,ChePin18}\footnote{While many estimators are introduced only for scalar-valued signals, 
they can also be used if $p>1$
by replacing  $\yaux(t)$ with
 $\yauxrvec^T \yaux(t)$
where $\yauxrvec\in \R^p$ is a fixed random vector.
\potential{
The randomness of $\yauxrvec$ guarantees the presence of all frequency components in $\yauxrvec^T \yaux(t)$ with probability $1$.
}
}.
Our controller requires 
frequency estimation
in the presence of the nonsmooth decaying part $y_0(t)$ of $\yaux(t)$ (i.e., the estimator is required to be input-to-state stable). 
One such adaptive estimator was introduced in~\cite{CarAst08} (see~\citel{CarAst08}{Rem.~3}).
The estimator in~\cite{CarAst08} is 
 compatible with our control scheme and it
only requires the minimal assumption of nonzero amplitudes.
During the past decade several estimators have been developed to improve the transient performance and robustness properties of earlier methods, e.g.,
in~\cite{CarGal16,ChePin18,VedVed21}.

Estimation of $(\gw_k)_{k=1}^{q}$  from $\yaux(t)$ requires that all frequencies 
appear in the non-decaying part of $\yaux(t)$.
This is generically true since the amplitudes corresponding to $\pm\gw_{k_0}$ are zero only if $\yrefamplk[k_0]$, $\wdistamplk[k_0]$, $\theta_{k_0}$ and $\varphi_{k_0}$ 
are related in a very specific way
through the identity
 $\Ptot(\pm i\gw_{k_0})c_e^{\pm k_0}=0$.%
\end{remark}%

\subsection{The Controller Tuning Algorithm}
\label{sec:TuningAlgorithm}

In this section we introduce an algorithm for
constructing 
 $(\hat{\gw}_k(\cdot))_{k=1}^q$
and 
 $K(\cdot)$ in the controller.
Even though several estimators
provide continuous-time estimates of $(\gw_k)_{k=1}^q$, we 
choose the estimates $(\hat{\gw}_k(\cdot))_{k=1}^q$ in the internal model to be piecewise constant functions which are updated at predefined time instances $0=t_0<t_1<t_2<\cdots$ (via sample-and-hold). 
This way we can guarantee stable closed-loop behaviour during the update intervals $[t_j,t_{j+1}]$
despite possible rapid changes 
in the frequency estimates.
The algorithm utilises the \emph{Estimate Admissibility Condition} defined below.

\begin{definition}
\label{def:FreqCond}
Let $\eps_f>0$ and $M_f>0$.
We say that
 $( \hat{\gw}_k)_{k=1}^q$ satisfy the \emph{Estimate Admissibility Condition} $\textup{EAC}(M_f,\eps_f)$ for the system $(A,B,C,D)$ if the following hold:
\begin{itemize}
\setlength{\itemsep}{.7ex}
\item $\eps_f\leq \abs{\hat{\gw}_k}\leq M_f$ for all $k\in \List{q}$.
\item $\abs{\hat{\gw}_k-\hat{\gw}_j}\geq \eps_f$ for all $k\neq j$.
\item $\abs{i\hat{\gw}_k\pm i\gl}\geq \eps_f$ 
 for every  
transmission zero $i\gl\in i\R$ of $(A,B,C,D)$ and
for all $k$.
\end{itemize}
\end{definition}

The algorithm uses
 Theorem~\ref{thm:Ae0Be0stabRLS} below to stabilize the pairs
$(A_s(t_j),B_s)$,
where
\eqn{
\label{eq:Ae0Be0defs}
  A_s(t)=\pmat{G_1(t)&G_2C\\0&A}, \quad B_s=\pmat{G_2D\\B}
}
with $G_1(t)
= \diag(0_p, \hat{\gw}_1(t)\Omega_p,\ldots, \hat{\gw}_q(t)\Omega_p)$.
The result guarantees that the stabilizing feedback gains $K^j$ are \emph{a priori} bounded and the stabilized semigroups satisfy uniform decay estimates with $M_s,\ga_s>0$ independent of $t_j$. 
Moreover, 
in this method the stabilizing gain
 of the infinite-dimen\-sional pair $(A,B)$ does not need to be recomputed when the frequencies in $G_1(t)$ are updated.
A similar method has been previously used in~\cite{HamPoh10,Pau16a} (see also \citel{Nat21}{Thm.~3.7}) for internal models with fixed frequencies.
The result uses notation
\eq{
H_0(i\gw) &= 
\frac{1}{2}
\pmat{ C_KR_K(i\gw) + C_KR_K(-i\gw)
\\ iC_KR_K(i\gw) -iC_KR_K(-i\gw)
}\\
B_0(i\gw) &= 
\frac{1}{2}
\pmat{
P_K(i\gw) + P_K(-i\gw)\\iP_K(i\gw) -iP_K(-i\gw)
},
}
where $R_K(\gl)=R(\gl,A+BK_{21})$,
 $C_K = \CL +DK_{21}$, and
 $P_K(\gl)=(\CL+DK_{21})R(\gl,A+BK_{21})B+D$.

\begin{theorem}
\label{thm:Ae0Be0stabRLS}
Let $\eps_f,M_f,r>0$ and let $Q_1\in\Lin(\C^{p(2q+1)})$ and $R_1\in \C^{m\times m}$ be positive definite.
Assume $K_{21}\in \Lin(X,\C^m)$ is 
such that $A+BK_{21} $ generates an exponentially stable semigroup $T_K(t)$
with growth bound $\gw_0(T_K(t))<0$.
 Assume $t_j\geq 0$ is such that $(\hat{\gw}_k(t_j))_{k=1}^q$ satisfy $\textup{EAC}(M_f,\eps_f)$ in Definition~\textup{\ref{def:FreqCond}} for the system $(A,B,C,D)$.
Define $H_j\in \Lin(X,\C^{p(2q+1)})$ and $B_{1j}\in \C^{p(2q+1)\times m}$ by 
\eq{
H_jx = 
 \pmat{ C_KR_K(0)x\\H_0(i\hat{\gw}_1(t_j))x \\
\vdots\\H_0(i\hat{\gw}_q(t_j))x} 
\quad \mbox{and}
\quad 
B_{1j} = 
 \pmat{ P_K(0)\\B_0(i\hat{\gw}_1(t_j)) \\
\vdots\\B_0(i\hat{\gw}_q(t_j))} 
}
Choose $K_1^j = -R_1\inv B_{1j}^\ast \Pi_{1j}\in \C^{m\times p(2q+1)}$ where $\Pi_{1j}\in \Lin(\C^{p(2q+1)})$ is the unique non-negative solution of 
\eq{
\MoveEqLeft[8] (rI+G_1(t_j))^\ast \Pi_{1j}+\Pi_{1j}(rI+G_1(t_j)) 
- \Pi_{1j}B_{1j}R_1\inv B_{1j}^\ast \Pi_{1j} = -Q_1.
}
If we choose $K^j = [K_1^j,K_{21}+K_1^jH_j]\in \Lin(Z,\C^m)$, then $\norm{K^j}\leq M_K$ for some $M_K>0$ independent of $t_j$.
Moreover,
 the semigroup $T_s^j(t)$ generated by $A_s(t_j)+B_sK^j$ 
 is exponentially stable so that 
for any $0<\ga_s<\min \set{r,-\gw_0(T_K(t))}$
 there exists $M_s,M_B>0$ (independent of $t_j$) such that
\eqn{
\label{eq:Tsunifbound}
\norm{T_s^j(t)}\leq M_s e^{-\ga_s t}, \qquad t\geq 0
}
and $\norm{R(\gl,A_s(t_j)+B_sK^j)B_s}\leq M_B$ for $\gl\in\C_+$.
If $(\hat{\gw}_k(t_j))_{k=1}^q$ satisfy $\EAC(M_f,\eps_f)$ for all $j$ large
and if $\max_k\abs{\hat{\gw}_k(t_j)- \gw_k}\to 0$ as $j\to\infty$, then
$\lim_{j\to\infty}\norm{K^j-K^\infty}=0 $,
 where $K^\infty =[K_1^\infty,K_{21}+K_1^\infty H_\infty]$ 
is obtained by replacing $(\hat{\gw}_k(t_j))_k$ with $(\gw_k)_k$ in $G_1(t_j)$, $B_{1j}$, and $H_j$.
\end{theorem}

The proof of Theorem~\ref{thm:Ae0Be0stabRLS} is presented in the Appendix.
In the tuning algorithm we denote by 
$0<\FE_1(t)<\cdots<\FE_q(t)$
the estimated frequencies 
computed based the signal $\yaux(t)$ by 
 the separate frequency estimator.
We suppose that Assumptions~\ref{ass:RLSstabilizability} and~\ref{ass:TZass} hold
and make the following additional assumptions on the parameters of the algorithm.

\begin{assumption}[Tuning Parameters]
  \label{ass:TuningPar}
~
\vspace{-.7ex}
\begin{itemize}
\item The sequence $0=t_0<t_1<\ldots$ of update times satisfies $\tau_1\leq t_j-t_{j-1}\leq \tau_2$ for some $\tau_1,\tau_2>0$ and all $j\in\N$.
\item The matrices $R_1\in \C^{m\times m}$ and $Q_1\in \Lin(\C^{p(2q+1)})$ are positive definite and $r>0$.
\item The frequency overlap parameter $\eps_f>0$ is suitably small and the upper bound $M_f>0$ for the frequencies is suitably large.
\item The operator $K_{21}\in \Lin(X,\C^m)$ 
is
 such that the semigroup generated by $A+BK_{21} $
is
 exponentially stable.
\item The initial frequency estimates $(\FE_k(0))_{k=1}^q$
satisfy $\EAC(M_f,\eps_f)$ for $(A,B,C,D)$.
\end{itemize}
\end{assumption}

The following Controller Tuning Algorithm constructs $G_1(\cdot)$
 (based on $(\hat{\gw}_k(\cdot))_{k=1}^q$)
and $K(\cdot)$ in the observer-based controller in Definition~\ref{def:ObsContr}.
The condition $\EAC(M_f,\eps_f)$ in Definition~\ref{def:FreqCond} is used in \textbf{Step~1} to detect if 
the frequency estimates nearly overlap or are close
 to the transmission zeros of $(A,B,C,D)$. In both cases the closed-loop stabilization becomes difficult, and therefore the the algorithm does not update the frequencies of the internal model if $\EAC(M_f,\eps_f)$ is violated (\textbf{Step~2} vs. \textbf{Step~3}).

\medskip

\noindent \textbf{The Controller Tuning Algorithm:}
\textit{
Choose $(t_j)_{j=0}^\infty$, 
$R_1$, $Q_1$, $\eps_f,M_f,r>0$, and $K_{21}$ as in Assumption~\textup{\ref{ass:TuningPar}}.
Set $j=0$.
}

\medskip

\noindent  \textbf{Step 1.} 
\textit{
Obtain
$0<\FE_1(t_j)<\cdots<\FE_q(t_j)$
 from the frequency estimator.
If 
$(\FE_k(t_j))_{k=1}^q$
 satisfy $\EAC(M_f,\eps_f)$ for $(A,B,C,D)$, then go to \textup{\textbf{Step~2}}. Otherwise go to \textup{\textbf{Step~3}}.
}

\medskip

\noindent  \textbf{Step 2 (Frequency update).} 
\textit{
Set $\hat{\gw}_k(t)\equiv \FE_k(t_j)$ 
  for $t\in[t_j,t_{j+1})$ and all $k$.
  Choose $K^j = [K_1^j,K_{21}+K_1^jH_j]\in \Lin(Z,\C^m) $ as in \textup{Theorem~\ref{thm:Ae0Be0stabRLS}} and
set $K(t)\equiv K^j$ for $t\in[t_j,t_{j+1})$.
  Increment $j$ to $j+1$ and go to \textup{\textbf{Step~1}}.
}

\medskip

\noindent  \textbf{Step 3 (No frequency update).} 
\textit{
Set $\hat{\gw}_k(t)\equiv \hat{\gw}_k(t_{j-1})$ 
and $K(t)\equiv K(t_{j-1})$ 
  for $t\in[t_j,t_{j+1})$ and all $k$.
Increment $j$ to $j+1$ and 
 go to \textup{\textbf{Step 1}}.
}

\medskip

Since the initial frequency estimates $(\FE_k(0))_{k=1}^q$ are assumed to 
satisfy  $\EAC(M_f,\eps_f)$, the tuning algorithm will proceed to \textbf{Step~2} when $j=0$, and therefore $G_1(\cdot)$ and $K(\cdot) $ are well-defined on $[0,\infty)$.
Our main result below shows that if the estimates $(\FE_k(t))_k$ converge to the true frequencies $(\gw_k)_k$ in~\eqref{eq:yrefwdist}, then
 the controller constructed with the above algorithm achieves output regulation of $\yref(t)$ and $\wdist(t)$.

\begin{theorem}
  \label{thm:ContrTuningThm}
Let  Assumptions~\textup{\ref{ass:RLSstabilizability}},~\textup{\ref{ass:TZass}}, and~\textup{\ref{ass:TuningPar}}  hold.
  Consider the controller $(\mc{G}_1(t),\mc{G}_2,K(t))$ in Definition~\textup{\ref{def:ObsContr}}, where $(\hat{\gw}_k(\cdot))_{k=1}^q$ and $K(\cdot)$ are based on the Controller Tuning Algorithm.
Assume
 that the true frequencies $(\gw_k)_{k=1}^q$ of $\yref(t)$ and $\wdist(t)$ 
satisfy $\EAC(\tilde{M}_f,\tilde{\eps}_f)$ for $(A,B,C,D)$ with some $\tilde{\eps}_f>\eps_f$ and $0<\tilde{M}_f<M_f$.

The controller satisfies $G_1(\cdot)\in \Lp[\infty](0,\infty;\Lin(Z_0))$ and $K(\cdot)\in \Lp[\infty](0,\infty;$ $\Lin(Z,\C^m))$.
If 
\potential{
$\wdist(t)$ and $\yref(t)$ and the intial conditions $x_0\in X$ and $z_0\in Z$ are such that%
}
$\abs{\FE_k(t)-\gw_k}\to 0$ as $t\to\infty$ for all $k$,
then
$U_e(t,s)$ is exponentially stable and
\eq{
\int_t^{t+1} \norm{y(s)-\yref(s)}^2ds\to 0 \qquad \mbox{as} \quad t\to\infty.
}
\end{theorem}

\begin{proof}
By construction $\abs{\hat{\gw}_k(t)}\leq M_f$ for all $k$, and thus $G_1(\cdot)\in \Lp[\infty](0,\infty;$ $\Lin(Z_0))$. 
Moreover, Theorem~\ref{thm:Ae0Be0stabRLS} shows that $K^j$ are uniformly bounded with respect to $j$, and thus $K(\cdot)\in \Lp[\infty](0,\infty;\Lin(Z,\C^m))$.
Assume now that $x_0\in X$, $z_0\in Z$, 
$\wdist(t)$, and $\yref(t)$ are such that 
$\abs{\FE_k(t)-\gw_k}\to 0$ as $t\to\infty$ for all $k$.
Since
 $\tilde{\eps}_f>\eps_f$ and $0<\tilde{M}_f<M_f$,
 there exists $N\in\N$ such that $(\FE_k(t_j))_{k=1}^q$ satisfy $\EAC(M_f,\eps_f)$ for all $j\geq N$. Thus the algorithm will go to \textbf{Step~2} for all $j\geq N$, and the frequency estimates of the internal model satisfy $\hat{\gw}_k(t)\to \gw_k$ as $t\to\infty$ for all $k$.
By Theorem~\ref{thm:Ae0Be0stabRLS} we have $K(t)\to K^\infty=[K_1^\infty,K_{21}+K_1^\infty H_\infty]\in \Lin(Z,\C^m)$ as $t\to\infty$, where $K_1^\infty$ and $H_\infty$ are obtained by replacing 
  $(\hat{\gw}_k(t_j))_k$ with $(\gw_k)_k$ in $G_1(t_j)$, $B_{1j}$, and $H_j$.
The claims will follow from Theorem~\ref{thm:ContrMain} provided that the semigroup generated by $A_s^\infty+B_sK^\infty$ is exponentially stable.
However, 
since $(\gw_k)_{k=1}^q$ satisfy $\EAC(M_f,\eps_f)$ by assumption, 
the stability of this semigroup follows directly from Theorem~\ref{thm:Ae0Be0stabRLS} when we replace
  $(\hat{\gw}_k(t_j))_k$ with $(\gw_k)_k$ in $G_1(t_j)$, $B_{1j}$, and $H_j$.
\end{proof}

The following lemma shows that for sufficiently large sampling intervals the evolution family $U_e(t,s)$ is always exponentially stable (independently of the behaviour of the frequency estimates  $(\FE_k(t))_k$).
Note that the required size of $\tau_1>0$ depends the other tuning parameters in Assumption~\ref{ass:TuningPar}.

\begin{lemma}
\label{lem:RobAnalUestab}
Let Assumptions~\textup{\ref{ass:RLSstabilizability}} and~\textup{\ref{ass:TuningPar}} hold.
There exists 
$\tau_1>0$ such that if $t_j-t_{j-1}\geq \tau_1$ for all $j\in\N$ in the Controller Tuning Algorithm,
 then there exist $M_e,\ga_e>0$ such that 
 $\norm{U_e(t,s)}\leq M_ee^{-\ga_e(t-s)}$ for all $t\geq s\geq 0$ and
 for any $\yref(t)$ and $\wdist(t)$ in~\eqref{eq:yrefwdist} and for any $x_0\in X$ and $z_0\in Z$.
\end{lemma}

\begin{proof}
Since $G_1(\cdot)$ and $K(\cdot)$ are piecewise constant we have $A_e(t)\equiv A_e(t_j)$ for all $t\in[t_j,t_{j+1})$ and $j\in\N_0$.
Thus if we denote by $T_e^j(t)$ the semigroup generated by $A_e(t_j)$, then for all $t\geq s\geq 0$ we have
 $U_e(t,s)=T_e^j(t-s)$  if $t,s\in [t_j,t_{j+1})$ for some $j\in\N_0$, and otherwise
\eq{
U_e(t,s) = 
T_e^j(t-t_j)T_e^{j-1}(t_j-t_{j-1})\cdots T_e^\ell(t_{\ell+1}-s)
}
where $j, \ell\in\N_0$ are such that $s\in [t_\ell, t_{\ell+1})$ and $t\in [t_j,t_{j+1})$.
Since $0<\tau_1\leq t_{j+1}-t_j\leq \tau_2$ for all $j\geq 0$ by assumption,
the evolution family $U_e(t,s)$ is exponentially stable
provided that there exists $M_{e0}\geq 0$ such that $\norm{T_e^j(t)}\leq M_{e0}$ for all $t\geq 0$ and $j\in\N_0$ and
 $\sup_{j\geq 0}\norm{T_e^j(t_{j+1}-t_j)}<1$.

 Theorem~\ref{thm:Ae0Be0stabRLS} and the Hille--Yosida theorem imply
 the existence of $M_s,\ga_s,$ $M_B,M_K>0$ such that
for $K^j$ in the Controller Tuning Algorithm 
we have 
 $\norm{R(\gl,A_s(t_j)+B_sK^j)}\leq M_s/(\re \gl+\ga_s)\leq M_s/\ga_s$, $\norm{R(\gl,A_s(t_j)+B_sK^j)B_s}\leq M_B$, and $\norm{K^j}\leq M_K$ for all $\gl\in\C_+$ and $j\in\N_0$
Using the similarity transform $Q_e$ in~\eqref{eq:Qetransform}
we have (similarly as in~\eqref{eq:QeAeform})%
\eq{
  Q_eA_e(t_j) Q_e\inv 
  &=\pmat{A_s(t_j)+B_sK^j&B_sK_2^j\\
  0&A+LC}.
}
The similarity, the triangular structure of $Q_eA_e(t_j)Q_e\inv$ and the norm estimates above imply
that
  there exists $M_R>0$ such that
$\sup_{\gl\in\C_+}\norm{R(\gl,A_e(t_j))}$ $\leq M_R$ for all $j\geq 0$.
Because of this, the Gearhart--Pr\"uss--Greiner theorem~\citel{EngNag00book}{Thm.~V.1.11} 
implies that
there exist $M_{e0},\ga_{e0}>0$ such that 
  $T_e^j(t)$ generated by $A_e(t_j)$  satisfy $\norm{T_e^j(t)}\leq M_{e0}e^{-\ga_{e0}t}$ for all $t\geq 0$ and $j\in \N_0$
(this uniform bound can be deduced, e.g., by applying~\citel{EngNag00book}{Thm.~V.1.11} to the semigroup $\diag(T_e^0(t),T_e^1(t),\ldots)$ on the Hilbert space $\lp[2](X_e)$).
This further implies that if we choose $\tau_1>0$ such that $M_{e0}e^{-\ga_{e0}\tau_1}<1$, then also 
$\norm{T_e^j(t_{j+1}-t_j)}\leq M_{e0}e^{-\ga_{e0}\tau_1}<1$ and the evolution family $U_e(t,s)$ is exponentially stable. Since $M_{e0}$ and $\ga_{e0}$ do not depend on $x_0$, $z_0$, $\wdist(t)$, and $\yref(t)$, we can choose
 $M_e,\ga_e>0$ as in the claim.
\end{proof}

\subsection{Robustness Analysis}
\label{sec:ContrRobAnalysis}

We conclude this section by analysing the robustness properties of the controller constructed in the Controller Tuning Algorithm.
The robustness properties depend on the chosen frequency estimation method --- especially on its capability of handling small persistent errors in $\yaux(t)$  --- but we can nevertheless present a general result for robustness analysis.
Throughout the section we assume that Assumptions~\ref{ass:RLSstabilizability}, \ref{ass:TZass} and \ref{ass:TuningPar} are satisfied.
We consider a perturbed regular linear system $(\tilde{A},\tilde{B},\tilde{C},\tilde{D})$ with parameters
\eqn{
\label{eq:PertSysParams}
\begin{cases}
\tilde{A}=A+\delta_A,\; 
\tilde{B}=B+\delta_B, \\
\tilde{C}=C+\delta_C, \;
\tilde{D}=D+\delta_D
\end{cases}
}
with $\delta_A\in \Lin(X)$, $\delta_B\in \Lin(\C^m,X)$, $\delta_C\in \Lin(X,\C^p)$ and $\delta_D\in \C^{p\times m}$.
We do not need to consider perturbations in $B_d$ and $D_d$ since these parameters were allowed to be unknown.
We begin by describing the effects of the perturbations on $\yaux(t)$.

\begin{lemma}
\label{lem:yauxpert}
Consider the perturbed system  $(\tilde{A},\tilde{B},\tilde{C},\tilde{D})$ in~\eqref{eq:PertSysParams} and 
 the controller in Definition~\textup{\ref{def:ObsContr}}.
Assume $K(\cdot)\in \Lp[\infty](0,\infty;\Lin(Z,\C^m))$ and $\hat{\gw}_k(\cdot)\in \Lp[\infty](0,\infty)$ for all $k$ are piecewise constant.
Then the auxiliary output $\yauxt(t)$ 
corresponding to the perturbed system satisfies
$\yauxt(t) = \yaux(t) + \ypert(t)$ for a.e. $t\geq 0$, where $\yaux(t)$ is as in Lemma~\textup{\ref{lem:yauxindep}} and 
\eqn{
\label{eq:ypertformula}
\ypert(t) = \CL\int_0^t T_L(t-s)
[\delta_{AC},\delta_{BD}K(s)]
x_e(s)ds
}
where $\delta_{AC}=\delta_A+L\delta_C$, $\delta_{BD}=\delta_B+L\delta_D$ and
 $x_e(t)$ is the state of the (perturbed) closed-loop system~\eqref{eq:CLS}.
\end{lemma}

\begin{proof}
Since $(\tilde{A},\tilde{B},\tilde{C},\tilde{D})$ is a regular linear system,  Theorem~\ref{thm:CLWellposedness} implies that the closed-loop system consisting of the perturbed system and the controller in Definition~\ref{def:ObsContr} has a well-defined mild state
$x_e(t)$.
Denote by $A_e(t)$ and $\tilde{A}_e(t)$ the closed-loop system operators corresponding to the nominal system $(A,B,C,D)$ and the perturbed system $(\tilde{A},\tilde{B},\tilde{C},\tilde{D})$, respectively.
Then $\Dom(\tilde{A}_e(t))=\Dom(A_e(t))$ and $\tilde{A}_e(t)=A_e(t) + \delta_e(t)$ for a.e. $t\geq0$ where 
\eq{
\gd_e(t)=\pmat{\gd_A&\gd_BK_1(t)&\gd_B K_2(t)\\
G_2 \gd_C&G_2\gd_DK_1(t)&G_2\gd_D K_2(t)\\
-L \gd_C&-L\gd_DK_1(t)&-L\gd_D K_2(t) }\in \Lin(X_e).
}
If we apply the similarity transform $Q_e$ in~\eqref{eq:Qetransform} to the perturbed closed-loop system~\eqref{eq:CLS}, we obtain
\eq{
\ddb{t}(Q_ex_e(t)) &= Q_eA_e(t)Q_e\inv (Q_ex_e(t))
+ Q_eB_e\wext(t) + Q_e\delta_e(t)x_e(t)
}
where
 $Q_e\gd_e(\cdot)x_e(\cdot)\in \Lploc[2](0,\infty;X_e)$ by Remark~\ref{rem:CLWPremark}.
The triangular structure of $Q_eA_e(t)Q_e\inv $ therefore implies that $\tilde{x}(t):=x(t)-z_2(t)$ is \emph{formally} a solution of 
\eqn{
\label{eq:yauxpertACP}
\hspace{-1ex}\dot{\tilde{x}}(t)
= A_L \tilde{x}(t)
 +B_{dL} \wext(t) + [\delta_{AC},\delta_{BD}K(t)]x_e(t)
}
with initial condition $\tilde{x}(0)=x_0-z_{20}$ and with $A_L=A+LC$. 
We will now prove that $\tilde{x}(t)$ is indeed a mild solution of~\eqref{eq:yauxpertACP}.
Denote by 
 $U_e(t,s)$ and $\tilde{U}_e(t,s)$ the evolution families in Theorem~\ref{thm:CLWellposedness} corresponding to the nominal and perturbed systems, respectively.
Let $0=t_0<t_1<t_2<\cdots$ be such that $A_e(t)\equiv A_e(t_j)$ and $\gd_e(t)\equiv \gd_e(t_j)$ for $t\in [t_j,t_{j+1})$. If we denote 
by $T_e^j(t)$ and $\tilde{T}_e^j(t)$ the semigroups generated by $A_e(t_j)$ and $\tilde{A}_e(t_j)$, respectively, then 
$U_e(t,s)$ and $\tilde{U}_e(t,s)$ are of the form given
 in the proof of Lemma~\ref{lem:RobAnalUestab}. Moreover, since $\tilde{A}_e(t_j)=A_e(t_j)+\gd_e(t_j)$ for all $j\in\N_0$, the perturbation formula in~\citel{EngNag00book}{Cor.~III.1.7} and a direct computation show that 
\eqn{
\label{eq:UeUepertRelationship}
\tilde{U}_e(t,s)x=U_e(t,s)x+\int_s^t U_e(t,r)\gd_e(r) \tilde{U}_e(r,s)xdr, 
}
for all $x\in X_e$ and $t\geq s\geq 0$.
Applying the similarity transformation $Q_e$ to~\eqref{eq:CLSstateoutput}  and~\eqref{eq:CLMapsDefPhi}
and using 
the relationship~\eqref{eq:UeUepertRelationship}
between $\tilde{U}_e(t,s)$ and $U_e(t,s)$ 
it is straightforward to confirm
 that $\tilde{x}(t)$ is the mild solution of~\eqref{eq:yauxpertACP}.
Since $\tilde{x}(t)\in \Dom(\CL)$ for a.e. $t\geq0$, analogous arguments as in the proof of Lemma~\ref{lem:yauxindep} show $\yauxt(t) = \CL \tilde{x}(t)+[D_d,\,-I]\wext(t)$. Comparing~\eqref{eq:yauxformACP} and \eqref{eq:yauxpertACP} shows that $\yauxt(t)=\yaux(t)+\ypert(t)$ for a.e. $t\geq 0$.
\end{proof}

Our main result below
shows that for sufficiently long sampling intervals in the Controller Tuning Algorithm
 the effect of small perturbations 
on $\yaux(t)$ will be small. Moreover, 
if the frequencies can be estimated with a sufficiently small asymptotic error, 
then the controller achieves output tracking \emph{in an approximate sense}, i.e., 
 with a small asymptotic error.

\begin{theorem}
\label{thm:ContrRobustnessMain}
Consider the perturbed system $(\tilde{A},\tilde{B},\tilde{C},\tilde{D})$ in~\eqref{eq:PertSysParams} and
let Assumptions~\textup{\ref{ass:RLSstabilizability}}, \textup{\ref{ass:TZass}} and~\textup{\ref{ass:TuningPar}} hold.
Let $\tau_1>0$ be as in Lemma~\textup{\ref{lem:RobAnalUestab}}
  and consider the controller
 in Definition~\textup{\ref{def:ObsContr}}, where $(\hat{\gw}_k(\cdot))_{k=1}^q$ and $K(\cdot)$ are based on the Controller Tuning Algorithm.

There exist $\eps_{stab},M_{pert}>0$ such that 
if 
\eqn{
\label{eq:ContrRobustnessDeltaSmallness}
c_\gd := \norm{\gd_A}+\norm{\gd_B}+\norm{\gd_C}+\norm{\gd_D}\leq \eps_{stab},
}
then for all $x_0\in X$, $z_0\in Z$,  $\wdist(t)$, and $\yref(t)$ 
we have
 $\yauxt(t)=\yaux(t)+\ypert(t)$,
where $\yaux(t)$ is as in Lemma~\textup{\ref{lem:yauxindep}} and 
\eq{
\MoveEqLeft\sup_{\tau\geq 0} 
\int_\tau^{\tau+1} 
\hspace{-1ex} \norm{\ypert(t)}^2dt
\leq  M_{pert}c_\gd^2 
( \norm{x_{e0}}^2+\norm{\wext(\cdot)}_\infty^2)
}
with $x_{e0}=[x_0,z_0]^T$ and $\wext(t)=[\wdist(t),\yref(t)]^T$.
If $C\in \Lin(X,\C^p)$, then
 $\norm{\ypert}_\infty^2\leq M_{pert}c_\gd^2( \norm{x_{e0}}^2+\norm{\wext(\cdot)}_\infty^2)$.

Assume
 $(\gw_k)_k$ in~\eqref{eq:yrefwdist} satisfy $\EAC(\tilde{M}_f,\tilde{\eps}_f)$ with some $\tilde{\eps}_f>\eps_f$ and $0<\tilde{M}_f<M_f$.
For any $\eps_{err}>0$ there exists $\gd_{err}>0$ such that if the perturbations satisfy~\eqref{eq:ContrRobustnessDeltaSmallness} and 
 $x_0\in X$, $z_0\in Z$, $\wdist(t)$, and $\yref(t)$ 
are such that
 $(\FE_k(t))_{k=1}^q$  satisfy
\eqn{
\label{eq:ContrRobustnessFEaccuracy}
\max_{k}\;\abs{\FE_k(t)-\gw_k}\leq \gd_{err}, \qquad \forall t\geq \tau_0
}
for some $\tau_0>0$, then 
\eq{
\limsup_{t\to\infty} \int_t^{t+1}\norm{y(s)-\yref(s)}^2ds\leq \eps_{err} \norm{w_e(\cdot)}_\infty^2.
}
\end{theorem}

\begin{proof}
Fix  $x_0\in X$, $z_0\in Z$, $\wdist(t)$, and $\yref(t)$ and
let $M_K>0$ be as in Theorem~\ref{thm:Ae0Be0stabRLS}.
Then
 $G_1(\cdot)$ and $K(\cdot)$  constructed in the Controller Tuning Algorithm 
 satisfy $\norm{K(\cdot)}_{\Lp[\infty]}\leq M_K$ and $\norm{G_1(\cdot)}_{\Lp[\infty]}\leq M_f$ and
thus
 by Theorem~\ref{thm:CLWellposedness}
 the closed-loop system corresponding to the perturbed system $(\tilde{A},\tilde{B},\tilde{C},\tilde{D})$  has a well-defined evolution family $\tilde{U}_e(t,s)$ and state $x_e(t)$.
We begin by introducing some notation.
We denote by $A_e(t)$ and $\tilde{A}_e(t)$ the closed-loop operators for
 $(A,B,C,D)$ and
 $(\tilde{A},\tilde{B},\tilde{C},\tilde{D})$, respectively.
Moreover, we denote by $T_e^j(t)$ and $\tilde{T}_e^j(t)$ the semigroups generated by $A_e(t_j)$ and $\tilde{A}_e(t_j)$, respectively.
For the proof of the last claim we
additionally assume that $(\gw_k)_k$ satisfy $\EAC(\tilde{M}_f,\tilde{\eps}_f)$. We then define  $\G_1^\infty$, $\G_{11}^\infty$ and $\Delta_{\G_1}(t)$ as in~\eqref{eq:ObsContrLimit} with
 $G_1^\infty = \diag(0_p,\gw_1\Omega_p,\ldots,\gw_q\Omega_p)$ 
 (i.e., $(\hat{\gw}_k(t))_k$ in $G_1(t)$ replaced with $(\gw_k)_k$) and define
$K^\infty=[K_1^\infty,K_{21}+K_1^\infty H_\infty]$ as in Theorem~\ref{thm:Ae0Be0stabRLS}.
Finally, we denote by $A_e^\infty$ and $\tilde{A}_e^\infty$ the operators in~\eqref{eq:CLlimitoperatorsAe} for
 $(A,B,C,D)$ and
 $(\tilde{A},\tilde{B},\tilde{C},\tilde{D})$, respectively, and denote the semigroups generated by these two operators with  $T_e(t)$ and $\tilde{T}_e(t)$, respectively.
Note that $\G_1^\infty$, $\G_{11}^\infty$, $\G_2$, and $K^\infty$ are independent of $x_0\in X$, $z_0\in Z$, $\yref(t)$, and $\wdist(t)$.

If $M_{e0},\ga_{e0}>0$ are as in the proof of Lemma~\ref{lem:RobAnalUestab},
we have $\norm{T_e^j(t)}\leq M_{e0}e^{-\ga_{e0}t}$ for all $t\geq 0$ and $j\in\N_0$ and
 the choice of $\tau_1>0$ implies
$M_{e0}e^{-\ga_{e0}\tau_1}<1$.
By construction $A_e^\infty$
 has the same structure
  as $A_e(t_j)$, 
 $j\in\N_0$, with $(\hat{\gw}_k(t_j))_k$ and $K^j$ replaced with $(\gw_k)_k$ and $K^\infty$, respectively, and  
$K^\infty$ in Theorem~\ref{thm:Ae0Be0stabRLS} is chosen similarly as $K^j$. Under the additional assumption that $(\gw_k)_k$ satisfy $\EAC(\tilde{M}_f,\tilde{\eps}_f)$, the frequencies 
$(\gw_k)_k$ satisfy
the assumptions of Theorem~\ref{thm:Ae0Be0stabRLS} 
with the same parameters as 
 $(\hat{\gw}_k(t_j))_k$, $j\in\N_0$ (both satisfy $\EAC(M_f,\eps_f)$ for $(A,B,C,D)$).  Therefore we can deduce as in the proof of Lemma~\ref{lem:RobAnalUestab} that also $\norm{T_e(t)}\leq M_{e0}e^{-\ga_{e0}t}$ for all $t\geq 0$.

We will now choose $\eps_{stab}>0$
so that~\eqref{eq:ContrRobustnessDeltaSmallness} implies the existence of
 $M_e,\ga_e>0$ (independent of $x_0$, $z_0$, $\wdist(t)$, $\yref(t)$)
such that $\norm{\tilde{U}_e(t,s)}\leq M_ee^{-\ga_e(t-s)}$ for $t\geq s\geq 0$
(our choice
 will also later
 guarantee the
 stability of
$\tilde{T}_e(t)$).
We have $\tilde{A}_e(t_j)=A_e(t_j)+\delta_e(t_j)$, $j\in\N_0$, where
$\delta_e(\cdot)\in \Lp[\infty](0,\infty;\Lin(X_e))$ 
is as in the proof of Lemma~\ref{lem:yauxpert}, and 
\eq{
 \norm{\delta_e(\cdot)}_{\Lp[\infty]}\leq
M_1\bigl( \norm{\gd_A}+\norm{\gd_B}+\norm{\gd_C}+\norm{\gd_D}\bigr)
 = M_1c_{\gd}
}
for some $M_1>0$ depending only on $L$ and $M_K$.
For a fixed $r_0\in (0,1)$ we choose $\eps_{stab}>0$ to be  small enough so that 
\eq{
\eps_{stab} \leq \frac{\ga_{e0}}{2M_{e0}M_1}
\quad \mbox{and} \quad 
M_{e0}e^{(-\ga_{e0}+M_{e0}M_1\eps_{stab})\tau_1}\leq r_0.
}
With this choice the condition
\eqref{eq:ContrRobustnessDeltaSmallness}
together with~\citel{EngNag00book}{Thm.~III.1.3} 
  implies that
$\norm{\tilde{T}_e^j(t)}\leq M_{e0}e^{-(\ga_{e0}/2)t}$ for all $t\geq 0$ and $j\in\N_0$
and 
$\sup_{j\geq 0}\norm{\tilde{T}_e^j(t_{j+1}-t_j)}\leq \sup_{j\geq 0}M_{e0}e^{(-\ga_{e0}+M_{e0}\norm{\gd_e(t_j)})(t_{j+1}-t_j)}\leq r_0<1$.
Since $\tilde{A}_e(t)\equiv\tilde{A}_e(t_j)$ for $t\in[t_j,t_{j+1})$ and $j\in\N_0$, 
the structure of $\tilde{U}_e(t,s)$ is analogous to that
 in the proof of Lemma~\ref{lem:RobAnalUestab}. 
Because the above estimates
for $\norm{\tilde{T}_e^j(t)}$ and $\norm{\tilde{T}_e(t_{j+1}-t_j)}$ 
 are uniform with respect to $(\gd_A,\gd_B,\gd_C,\gd_D)$ satisfying~\eqref{eq:ContrRobustnessDeltaSmallness}, 
we have (similarly as
 in the proof of Lemma~\ref{lem:RobAnalUestab}) that there exist $M_e,\ga_e>0$ such that $\norm{\tilde{U}_e(t,s)}\leq M_ee^{-\ga_e(t-s)}$ for all $t\geq s\geq 0$ and for any $(\gd_A,\gd_B,\gd_C,\gd_D)$ for which~\eqref{eq:ContrRobustnessDeltaSmallness} holds.

We will now prove the claims concerning $\ypert(t)$.
If $M_e,\ga_e>0$ are as above,
\potential{
\eqref{eq:CLSstateoutputstate},~\eqref{eq:CLMapsDefPhi} and $B_e(t)\equiv B_e\in \Lin(\C^{n_d+p},X_e)$ imply that for all $t\geq 0$ 
\eq{
\norm{x_e(t)}
&=\norm{\tilde{U}_e(t,0)x_{e0} + \int_0^t\tilde{U}_e(t,s)B_e(s)\wext(s)ds}\\
&\leq M_e \norm{x_{e0}} +  \frac{M_e \norm{B_e}}{\ga_e} \norm{w_e(\cdot)}_{\infty}.
}
Thus%
}
 $\norm{x_e(\cdot)}_{\infty}\leq M_2 (\norm{x_{e0}}+\norm{w_e(\cdot)}_\infty)$ for a constant $M_2\geq 0$ independent of $x_0$, $z_0$, $\wdist(t)$ and $\yref(t)$.
 Lemma~\ref{lem:yauxpert} shows that  $\ypert(t)= (\F_L u)(t)$ for a.e. $t\geq 0$,
where $\F_L$ is the input-output map of the exponentially stable regular linear system $(A+LC,I,C,0)$ and where $u(\cdot)=[\gd_{AC},\gd_{BD}K(\cdot)]x_e(\cdot)\in\Lploc[2](0,\infty;X)$ by Remark~\ref{rem:CLWPremark}.
Since $\norm{u}_{\Lp[2](\tau,\tau+1)}
\leq (M_K+1)\max\set{\norm{\gd_{AC}},\norm{\gd_{BD}}}\norm{x_e(\cdot)}_{\infty}\leq (M_K+1)M_2M_3c_\gd (\norm{x_{e0}}+\norm{w_e(\cdot)}_\infty)$ for all $\tau\geq 0$
and for a constant $M_3\geq 0$ depending only on $\norm{L}$,
 the first estimate for $\norm{\ypert(\cdot)}$ for some $M_{pert}\geq 0$ follows from
 Lemma~\ref{lem:IOmapconvproperties}(a). 
If $C\in \Lin(X,\C^p)$ and if $M_L,\ga_L>0$ are such that $\norm{T_L(t)}\leq M_Le^{-\ga_Lt}$ for all $t\geq 0$, then 
 the second claim follows from~\eqref{eq:ypertformula} and
a direct estimate 
\eq{
 \norm{\ypert(t)} 
& \hspace{-.3ex} \leq \hspace{-.3ex} \norm{C} \hspace{-.2ex} \max\set{\norm{\gd_{AC}}, \hspace{-.2ex}\norm{\gd_{BD}}}
\frac{M_L(M_K \hspace{-.25ex}+ \hspace{-.25ex}1)}{\ga_L}\norm{x_e(\cdot)}_{\infty}.
}

To prove the last claim we will apply Theorem~\ref{thm:ORPmainNonconv} to the 
perturbed
 system $(\tilde{A},\tilde{B},\tilde{C},\tilde{D})$.
Assume $(\gw_k)_k$ satisfy $\EAC(\tilde{M}_f,\tilde{\eps}_f)$ and let $\eps_{stab}>0$ be as above.
 It is easy to see that the perturbations~\eqref{eq:PertSysParams}
lead to bounded perturbations of
$(A_e^\infty,B_{ee}^\infty,C_{ee}^\infty,D_{ee}^\infty)$ in~\eqref{eq:CLlimitoperators} 
with norm bounds depending on $c_\gd$, $M_K$ and $\norm{L}$.
Moreover, by the choice of $\eps_{stab}$ the perturbation in $A_e^\infty$ has norm at most $\ga_{e0}/(2M_{e0})$ when~\eqref{eq:ContrRobustnessDeltaSmallness} holds.
Thus Remark~\ref{rem:ORPnonconvDeltaDependence} and Lemma~\ref{lem:IOmapperts} imply that 
$\Merr$ and $\gd_0$ in Theorem~\ref{thm:ORPmainNonconv} can be chosen to hold for
 all perturbations $(\gd_A,\gd_B,\gd_C,\gd_D)$ satisfying~\eqref{eq:ContrRobustnessDeltaSmallness}.
The controller $(\G_1(t),\G_2(t),K(t))$ satisfies the assumptions of Theorem~\ref{thm:ORPmainNonconv} since Assumption~\ref{ass:ContrAss} is satisfied by construction, $(\G_1^\infty,\G_2)$ satisfy~\eqref{eq:Gconds} in Remark~\ref{rem:Gconds}
 (similarly as in the proof of Theorem~\ref{thm:ContrMain}), and as shown above, both $\tilde{U}_e(t,s)$ and $\tilde{T}_e(t)$ are 
 exponentially stable whenever~\eqref{eq:ContrRobustnessDeltaSmallness} holds.
The definition of $\Delta_{\G}(t)$ in~\eqref{eq:ObsContrLimit}  and $\norm{G_1(t)-G_1^\infty}=\max_{k}\abs{\hat{\gw}_k(t)-\gw_k}$ imply
that $\gd_{\G}(t)$ in Theorem~\ref{thm:ORPmainNonconv} satisfies
\eq{
\gd_{\mc{G}}(t)
 &= \max \set{\norm{\Delta_{\mc{G}_1}(t)},\norm{K(t)-K^\infty}}
 = \norm{\Delta_{\mc{G}_1}(t)}
\\
&\leq  \norm{K^j-K^\infty} + \max_k\, \abs{\hat{\gw}_k(t)-\gw_k} ,
}
where $j\in\N_0$ is such that 
$t\in [t_j,t_{j+1})$. 
Since $(\gw_k)_k$ satisfy $\EAC(\tilde{M}_f,\tilde{\eps}_f)$ with $\tilde{M}_f<M_f$ and $\tilde{\eps}_f>\eps_f$, by choosing a sufficiently small $\gd_{err}>0$ we can guarantee that 
if~\eqref{eq:ContrRobustnessFEaccuracy} holds, then
$(\mu_k(t_j))_k$ satisfy $\EAC(M_f,\eps_f)$ for all $j\in\N$ such that $t_j\geq \tau_0$. Such choice guarantees that  $(\hat{\gw}_k(t_j))_k$ and $K(t_j)$ are updated whenever $t_j\geq \tau_0$, and thus we also have $\max_k \abs{\hat{\gw}_k(t)-\gw_k}\leq \gd_{err}$ for all $t\geq \tau_0+\tau_2$.
Theorem~\ref{thm:Ae0Be0stabRLS} shows that 
$\lim_{j\to\infty}K^j= K^\infty$
if $\lim_{j\to\infty}\max_k \abs{\hat{\gw}_k(t)-\gw_k}= 0$ 
and thus
  $\gd_{\mc{G}}(t)$ for $t\in[t_j,t_{j+1})$ can be made arbitrarily small by requiring that $\max_k\abs{\hat{\gw}_k(t_j)-\gw_k}$ is small.
For any $\eps_{err}>0$ we can now combine the above properties to choose 
 $\gd_{err}>0$ (independent of $x_0$, $z_0$,  $\wdist(t)$, and $\yref(t)$) such that if~\eqref{eq:ContrRobustnessFEaccuracy} holds for some $\tau_0>0$, then $\limsup_{t\to\infty}\norm{\gd_{\G}(\cdot)}_{\Lp[\infty](t,\infty)}^2<\max\set{\gd_0^2,\eps_{err}/\Merr}$. 
We then have from Theorem~\ref{thm:ORPmainNonconv} 
that~\eqref{eq:ContrRobustnessDeltaSmallness} and~\eqref{eq:ContrRobustnessFEaccuracy} for some $\tau_0>0$ imply
\ieq{
\limsup_{t\to\infty}\, \norm{e(\cdot)}_{\Lp[2](t,t+1)}^2
\leq \Merr \norm{w_e(\cdot)}_\infty^2 \limsup_{t\to\infty}\norm{\gd_{\G}(\cdot)}_{\Lp[\infty](t,\infty)}^2
\leq \eps_{err} \norm{w_e(\cdot)}_\infty^2.
}
\end{proof}

\section{Adaptive Regulation for a Heat Equation}
\label{sec:simulation}

In this example we study the output regulation problem for 
 a one-dimen\-sional boundary controlled reaction-diffusion equation. The system on $\xi\in(0,1)$ has the form 
\eq{
\pd{v}{t}(\xi,t) &= \pd[2]{v}{\xi}(\xi,t) + \gg(\xi)v(\xi,t) + b_d(\xi) \wdistk[1](t)\\
-\pd{v}{\xi}(0,t)&=u(t)+\wdistk[2](t), \quad \pd{v}{\xi}(1,t)=\wdistk[3](t)\\
y(t)&=v(1,t), \qquad 
v(\xi,0)=v_0(\xi), 
}
where $v(\xi,t)$ describes temperature at the point $\xi\in (0,1)$ and at time $t>0$.
 The boundary input $u(t)$ acts at $\xi=0$, the output $y(t)$ is the temperature measurement at $\xi=1$, and $\wdistk[2](t)$ and $\wdistk[3](t)$ are boundary disturbances.
The reaction term with profile $\gg(\cdot)$ is unmodelled and we consider it as perturbation in the system.
The disturbance input profile $b_d(\xi)$ is unknown.
The system defines a regular linear system on $X=\Lp[2](0,1)$ with state $x(t)=v(\cdot,t)$. The full disturbance input is defined as $\wdist(t)=[\wdistk[1](t),\wdistk[2](t),\wdistk[3](t)]^T\in\R^3$. Since the boundary disturbances $\wdistk[2](t)$ and $\wdistk[3](t)$ are smooth functions, we can apply a change of variables as
 in~\citel{CurZwa20book}{Sec.~10.1,~Ex.~10.1.7}
 to express the heat equation as a regular linear system with bounded $B_d$, $D_d\in \R^{1\times 3}$ and a modified initial state.

We construct a controller for output regulation of 
$\yref(t)=  0.2\sin(0.5t+0.5)+0.4\sin(6t+0.5)$ and  $\wdist(t) = [\cos(1.5t+0.5),\sin(0.5t+0.2),\cos(1.5t-0.4)]^T$, both assumed to be unknown. 
In the simulation we consider $\gg(\xi)=1.5\sin(0.5\pi \xi)$ and $b_d(\xi)=\cos(3\xi)$ (these are not used in the controller design).
We assume the number of nonzero frequencies $q=3$ is known.
The stabilizing parameters are chosen as $K_{21}x = 2\int_0^1x(\xi)d\xi$ for $x\in \Lp[2](0,1)$, $L\equiv -4\in \Lp[2](0,1)$.
The system does not have transmission zeros on $i\R$.
We use $\eps_f=0.2$ and $M_f=30$ in the Estimate Admissibility Condition, and $r=0.2$, $R=1\in\R$, and $Q=I\in\R^{3\times 3}$.
The frequencies will also not be updated if $(\FE_k(t_j))_{k=1}^3$ are complex or negative.

\begin{figure}[h!]
\begin{center}
\includegraphics[width=0.65\linewidth]{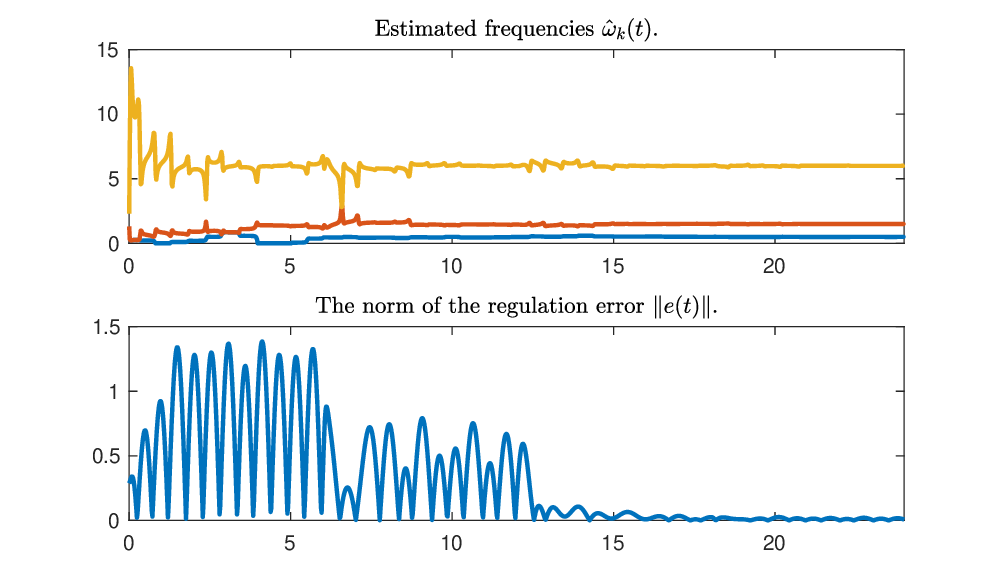}
\end{center}
\caption{Controlled heat equation with the estimator in~\cite{CarAst08}.}
\label{fig:heatmain}
\end{figure}

We use the adaptive estimator from~\cite{CarAst08} 
with parameters
 ``$\gamma_1=0.005$'', ``$\gamma_2=10$'', and ``$\set{k_i}_i$'' being the coefficients of the Hurwitz polynomial $(\lambda+2)^{2\cdot 3-1}$.  Initial frequency estimates are chosen as $\FE_k(0)=k\in\R$ for $k\in \List{3}$.
 The simulations are implemented using Finite Difference with $100$ points on $[0,1]$.
The input-to-state stability of the estimator and Theorem~\ref{thm:ContrRobustnessMain} imply that for sufficiently long update intervals the effect of the reaction term with a small $\norm{\gg}_{\Lp[2]}$ will lead to approximate output tracking with a small asymptotic error.
Figure~\ref{fig:heatmain} shows the behaviour of the frequency estimates and the regulation error 
for the update sequence $t_j=6j$, $j\in\N_0$, and initial states $x_0(\cdot,0)\equiv 0$ and $z_0=0\in Z$.

\subsection*{Acknowledgement}

The authors would like to thank Dr. Petteri Laakkonen for advice on transmission zeros and Prof. Dr. Andre Ran
 for helpful comments 
 related to Lemma~\ref{lem:IMstabRiccbound}.

\appendix

\section{Appendix}

\begin{lemma}
\label{lem:IOmapconvproperties}
Assume that $(U,\Phi,\Psi,\F)$ is a well-posed nonautonomous system in the sense of~\textup{\citel{Sch02}{Def.~3.6}}. If 
 there exist $M,\gw>0$ such that $\norm{U(t,s)}\leq M e^{-\gw (t-s)}$ for all $t\geq s\geq 0$, then for $u\in \Lploc[2](0,\infty;U)$
the following hold.
\begin{itemize}
\item[\textup{(a)}] 
There exists $M_0>0$ (independent of $u$) such that 
$\displaystyle\sup_{\tau\geq 0} \norm{\F u}_{\Lp[2](\tau,\tau+1)} \leq M_0\sup_{\tau\geq 0}\norm{u}_{\Lp[2](\tau,\tau+1)}$.
\item[\textup{(b)}]
There exists $M_1>0$ (independent of $u$) such that 
\eq{
\limsup_{\tau\to \infty}\;  \norm{\F u}_{\Lp[2](\tau,\tau+1)} \leq M_1\limsup_{\tau\to\infty}\;\norm{u}_{\Lp[2](\tau,\tau+1)}.
}
\item[\textup{(c)}]
If $\displaystyle \lim_{\tau\to \infty}\norm{ u}_{\Lp[2](\tau,\tau+1)} =0$,
  then 
 $\displaystyle \lim_{\tau\to \infty}\norm{\F u}_{\Lp[2](\tau,\tau+1)}= 0$.

\item[\textup{(d)}]
If
 $\sup_{\tau\geq 0} e^{\ga \tau}\norm{u}_{\Lp[2](\tau,\tau+1)}<\infty$ 
for some $0<\ga<\gw$, then  
$t\mapsto e^{\gb t}\norm{(\F u)(t)}\in\Lp[2](0,\infty)$
 for any $0<\gb<\ga$.
\end{itemize}
\end{lemma}

\begin{proof} 
Let $u\in \Lploc[2](0,\infty;U)$.
It is clearly sufficient to prove the claims when $\tau$ is replaced by $n\in\N_0$.
The estimates in~\citel{Sch02}{Lem.~3.7} show
 that there exists a constant $M_2>0$ (depending only on $(U,\Phi,\Psi,\F)$) such that 
\eq{
\norm{\F u}_{\Lp[2](n,n+1)} \leq M_2 (a\ast b)_n, \qquad \forall n\in\N_0,
}
where $(a\ast b)_n$ is the $n$th element 
of the convolution of 
$a=(a_k)_{k=0}^\infty$ with $a_k =e^{-\gw k}$ and 
 $b=(b_k)_{k=0}^\infty$ with $b_k = \norm{u}_{\Lp[2](k,k+1)}$.
The constant $M_2>0$ is determined by $M$ and $\gw$ and the uniform (w.r.t $s\geq 0$) bounds for 
$\norm{\Phi_{\cdot,s}}_{\Lin(\Lp[2](s,s+1),X)}$, $\norm{\Psi_s}_{\Lin(X,\Lp[2](s,s+1))}$, and $\norm{\F_s}_{\Lin(\Lp[2](s,s+1),\Lp[2](s,s+1))}$.

Since $a\in \lp[1](\R)$ and
 $\norm{b}_{\lp[\infty]}= \sup_{k\geq 0}\norm{u}_{\Lp[2](k,k+1)}$,
the claim in part (a) holds since
$\norm{a\ast b}_{\lp[\infty]}\leq \norm{a}_{\lp[1]}\norm{b}_{\lp[\infty]}$
by the Young's inequality for convolutions.
Part (c) follows from~(b).
To prove (b)
we assume $\limsup_{\tau\to\infty} \norm{u}_{\Lp[2](\tau,\tau+1)}$ $<\infty$ (otherwise the claim is trivial). Since $u\in \Lploc[2](0,\infty;U)$ we have $\norm{b}_{\lp[\infty]}=\sup_{k\geq 0}\norm{u}_{\Lp[2](k,k+1)}<\infty$.
If $n,n_0\in\N$ satisfy $n_0<n$,
then $a_{n-k}=e^{-\gw(n-n_0+1)}a_{n_0-1-k}$ and
\eq{
(a\ast b)_n 
&=e^{-\gw(n-n_0+1)}(a\ast b)_{n_0-1} + (a\ast \tilde{b})_{n-n_0},
}
where $\tilde{b} = (b_{k+n_0})_{k=0}^\infty$. 
Young's inequality thus implies
\eq{
\abs{(a\ast b)_n} 
&\leq e^{-\gw(n-n_0+1)}\norm{a\ast b}_{\lp[\infty]} + \norm{a\ast \tilde{b}}_{\lp[\infty]}\\
&\leq e^{-\gw(n-n_0+1)}\norm{a}_{\lp[1]}\norm{ b}_{\lp[\infty]} + \norm{a}_{\lp[1]} \norm{\tilde{b}}_{\lp[\infty]}.
}
If we choose $n_0=\floor{n/2}$, then
the properties of the limit supremum and
 $\norm{\tilde{b}}_{\lp[\infty]}= \sup_{k\geq \floor{n/2}}\norm{u}_{\Lp[2](k,k+1)}$ 
 imply (b).

Finally, to prove (d)
we note that
\ieq{
e^{\ga n}\norm{\F u}_{\Lp[2](n,n+1)} \leq M_2 (a_\ga\ast b_\ga)_n
}
for all $n\in\N_0$,
where 
$a_\ga=(e^{-(\gw-\ga)k})_{k=0}^\infty\subset \R$ and 
 $b_\ga=(e^{\ga k}\norm{u}_{\Lp[2](k,k+1)})_{k=0}^\infty\subset \R$. Since $0<\ga<\gw$, we have $a_\ga\in\lp[1](\R)$ and our assumptions imply $b_\ga\in \lp[\infty](\R)$. Thus
Young's inequality 
 implies 
$\norm{a_\ga\ast b_\ga}_{\lp[\infty]}\leq 
 \norm{a_\ga}_{\lp[1]} \norm{b_\ga}_{\lp[\infty]}$
and we have $\sup_{n\geq 0}e^{\ga n}\norm{\F u}_{\Lp[2](n,n+1)}<\infty$. This implies the claim for any $0<\gb<\ga$.
\end{proof}

\begin{lemma}
\label{lem:IOmapperts}
Let $(A,B,C,D)$ be a regular linear system. 
Assume that there exist $M,\ga>0$ such that the semigroup $T(t)$ generated by $A$ satisfies $\norm{T(t)}\leq Me^{-\ga t}$ for $t\geq 0$.
If we denote the extended input, output and input--output maps of a perturbed system $(\tilde{A},\tilde{B},\tilde{C},\tilde{D})$ by $\tilde{\Phi}$, $\tilde{\Psi}$, and $\tilde{\F}$, respectively, then for any $\eps\in (0,\ga/M)$ and $\kappa>0$  we have
\eq{
\sup_{(\tilde{A},\tilde{B},\tilde{C},\tilde{D})\in\Omega(\eps,\kappa)}\left( \norm{\tilde{\Phi}}+\norm{\tilde{\Psi}}+\norm{\tilde{\F}} \right)<\infty,
}
where
 $\Omega(\eps,\kappa) = \setm{(A+\gd_A,B+\gd_B,C+\gd_C,D+\gd_D)}{\norm{\gd_A}\leq \eps ~\mbox{and}~ \norm{\gd_B}+\norm{\gd_C}+\norm{\gd_D}\leq \kappa}$.
\end{lemma}

\begin{proof}
Denote 
the extended input, output and input--output maps of $(A,B,$ $C,D)$ by $\Phi$, $\Psi$, and $\F$, respectively. Let $\eps\in (0,\ga/M)$ be fixed.
We begin by considering perturbations in operator $A$ only, i.e.,
$(\tilde{A},B,C,D) =(A+\gd_A,B,C,D)$, where $\norm{\gd_A}\leq \eps$.
Denote the input, output, and input--output maps
of the extended system $(A,[B,I],\pmatsmall{C\\I},\pmatsmall{D&0\\0&0})$ by
\eq{
\Phi_e=[\Phi,\Phi_I], \quad \Psi_e=\pmat{\Psi\\\Psi_I}, \quad \F_e=\pmat{\F&\F_{CI}\\\F_{IB}&\F_{II}}.
}
By the results in~\citel{Wei94}{Sec.~7}, applying
 an admissible output feedback
 $u_e(t)=\Delta y_e(t)+\tilde{u}_e(t)$ 
 with $\Delta=\pmatsmall{0&0\\0&\gd_A}$
leads to the regular linear system
$(A+\gd_A,[B,I],\pmatsmall{C\\I},\pmatsmall{D&0\\0&0})$ with input map $\bar{\Phi}_e$, output map $\bar{\Psi}_e$, and input-output $\bar{\F}_e$.
We have from~\citel{Wei94}{Thm.~6.1} that
$\bar{\F}_e = (I-\F_e \Delta)\inv \F_e$.
Since $(\F_{II}u)(t)=\int_0^t T(t-s)u(s)ds$, 
\citel{AreBat11book}{Prop.~1.3.5(a)} implies 
$\norm{\F_{II}}_{\Lin(\Lp[2](0,\infty))}\leq \norm{T(\cdot)}_{\Lp[1](0,\infty)}\leq M/\ga$.
Thus $\norm{\F_{II}\gd_A}\leq M\norm{\gd_A}/\ga\leq M\eps/\ga<1$ and
\eq{
(I-\F_e \Delta)\inv 
&= \pmat{I&-\F_{CI}\gd_A\\0&I-\F_{II}\gd_A}\inv 
= \pmat{I&\F_{CI}\gd_AQ\\0&Q} 
}
where $Q=(I-\F_{II}\gd_A)\inv$.
This implies
that we have
 $\norm{(I-\F_e\Delta)\inv}
\leq (1+(1+\norm{\F_{CI}\gd_A}^2)\norm{Q}^2)^{1/2}
$.
Finally, the estimate 
$\norm{Q}\leq \ga/(\ga-M\eps)$ and the formulas 
$\bar{\Phi}_e = \Phi_e(I+\Delta \bar{\F}_e)$ 
and
$\bar{\Psi}_e = (I+ \bar{\F}_e\Delta)\Psi_e$ in~\citel{Wei94}{Rem.~6.5} imply that
\eqn{
\label{eq:IOmappertBar}
\sup_{\norm{\gd_A}\leq \eps}\Bigl( \norm{\bar{\Phi}_e} + \norm{\bar{\Psi}_e} + \norm{\bar{\F}_e} \Bigr)<\infty.
}

We will now consider perturbed systems satisfying
$(\tilde{A},\tilde{B},\tilde{C},\tilde{D}) =(A+\gd_A,B+\gd_B,C+\gd_C,D+\gd_D)\in \Omega(\eps,\kappa)$ with $\kappa>0$.
It is easy to verify that
$\tilde{\Psi}x = [I,\gd_C] \bar{\Psi}_e x$, 
\eq{
\tilde{\Phi}u = \bar{\Phi}_e \pmat{u\\\gd_B u} 
\quad \mbox{and} \quad
\tilde{\F}u = \pmat{I&\gd_C} \bar{\F}_e\pmat{u\\\gd_B u} + \gd_D u.
}
Since $\norm{\gd_A}\leq \eps$ and $\norm{\gd_B}+\norm{\gd_C}+\norm{\gd_D}\leq \kappa$, the claim follows directly from~\eqref{eq:IOmappertBar}.
\end{proof}

The following 
corollary of
the continuity of the solutions of Riccati equations is essential for the proof of Theorem~\ref{thm:Ae0Be0stabRLS}. 
To the best of the authors' knowledge, this result is new. 

\begin{lemma}
\label{lem:IMstabRiccbound}
Let $r>0$, $Q\in\C^{n\times n}$ and $R\in \C^{m\times m}$ satisfy $Q>0$ and $R>0$. Let $\Omega\subset \R^q$ be a compact set and 
let
$\gd\mapsto A_\gd:\Omega\to \C^{n\times n}$ and $\gd\mapsto B_\gd:\Omega\to \C^{n\times m}$ be continuous functions such that the pair $(A_\gd,B_\gd)$ is controllable for all $\gd\in\Omega$.
If we define $K_\gd=-R\inv B_\gd^\ast\Pi_\gd$, $\gd\in\Omega$, 
 where $\Pi_\gd\in \C^{n\times n}$ are the unique non-negative solutions of
\eq{
(rI+A_\gd)^\ast \Pi_\gd + \Pi_\gd (rI+A_\gd) - \Pi_\gd B_\gd R\inv B_\gd^\ast \Pi_\gd =-Q,
}
then there exist $M,M_K>0$ such that $\norm{K_\gd}\leq M_K$ and 
$\norm{e^{(A_\gd+B_\gd K_\gd)t}}\leq Me^{-rt}$ for all $t\geq 0$ and $\gd\in\Omega$.
\end{lemma}

\begin{proof}
The claim is trivially true if $\Omega$ is empty.
Let $\gd\in\Omega$.
The assumptions imply that
$\Pi_\gd$ exists and is unique, and
 for all $x\in \C^n$ we have
\ieq{
 2\re \iprod{(rI+A_\gd+B_\gd K_\gd )x}{\Pi_\gd x}
 =\iprod{- \Pi_\gd B_\gd R\inv B_\gd^\ast \Pi_\gd x  -Q x}{x}
\leq 0.
}
Moreover, under our assumptions $\Pi_\gd$ is positive definite.
Therefore $rI + A_\gd+B_\gd K_\gd $ is dissipative with respect to the inner product
 $\iprod{\cdot}{\cdot}_\gd := \iprod{\cdot}{\Pi_\gd \cdot}_{\C^n}$ on $\C^n$.
 Thus if we define $\norm{x}_\gd :=\norm{\Pi_\gd^{1/2}x}_{\C^n}$  for $x\in \C^n$, then 
$\norm{e^{(rI+A_\gd+B_\gd K_\gd )t}x}_\gd\leq \norm{x}_\gd$ for all $t\geq 0$.
\potential{The definition of $\norm{\cdot}_\gd$ now implies that} 
$\norm{e^{(A_\gd+B_\gd K_\gd )t}}\leq \norm{\Pi_\gd^{1/2}}\norm{\Pi_\gd^{-1/2}}e^{-rt}$ for all $t\geq 0$.

Our aim is 
 to show that $\sup_{\gd\in\Omega}\norm{\Pi_\gd^{1/2}}\norm{\Pi_\gd^{-1/2}}<\infty$. 
By~\citel{Sun98}{Thm.~3.1}
 the non-negative matrix $\Pi_\gd$ is a continuous function of the matrices $A_\gd$ and $B_\gd$ when $\gd$ is restricted to $\Omega$. Therefore the function $\gd\mapsto \Pi_\gd$ is continuous on $\Omega$, and since $\Pi_\gd$ and $\Pi_\gd^{1/2}$ are nonsingular for all $\gd\in\Omega$, also $\gd\mapsto \Pi_\gd^{1/2}$ and $\gd\mapsto \Pi_\gd^{-1/2}$ are continuous on $\Omega$. Since $\Omega$ is compact, these functions are uniformly continuous and $\norm{\Pi_\gd^{1/2}}$ and $\norm{\Pi_\gd^{-1/2}}$ are uniformly bounded with respect to $\gd\in\Omega$. 
Thus the claims hold with
 $M_K:=\norm{R\inv}\max_{\gd\in\Omega}\norm{B_\gd}\norm{\Pi_\gd}<\infty$ and $M:=\max_{\gd\in\Omega}\norm{\Pi_\gd^{1/2}}\norm{\Pi_\gd^{-1/2}}<\infty$.
\end{proof}

\begin{proof}[Proof of Theorem~\textup{\ref{thm:Ae0Be0stabRLS}}]
Since $K_{21}\in \Lin(X,\C^m)$,
 $(A+BK_{21},B,C+DK_{21},D)$ is an exponentially stable regular linear system and $\gl\mapsto (\CL+DK_{21})R(\gl,A+BK_{21})$ and $\gl\mapsto P_K(\gl)$ are continuous functions on $[-iM_f,iM_f]\subset i\R$. Thus
 $\norm{H_j}\leq M_H$ and $\norm{B_{1j}}\leq M_{B1}$
for some $M_H,M_{B1}>0$ independent of $t_j$.
Definition~\ref{def:TransmissionZero} implies that the transmission zeros of $(A,B,C,D)$ 
 are zeros of $\gl\mapsto \det(P_K(\gl)P_K(\gl)^\ast)$, which is analytic on $\setm{\gl\in\C}{\re\gl>\gw_0(T_K(t))}$. Thus the transmission zeros of $(A,B,C,D)$ on $i\R$ are a (possibly empty) discrete set with no finite accumulation points.

By assumption,
$(\hat{\gw}_k(t_j))_{k=1}^q$ satisfy $\textup{EAC}(M_f,\eps_f)$. Therefore 
the set $ \set{0}\cup \set{\pm i\hat{\gw}_k(t_j)}_{k=1}^q$ does not contain any transmission zeros of $(A,B,C,D)$ and  $P_K(\pm i\hat{\gw}_k(t_j))^\ast$ and $P_K(0)^\ast$ are injective.
Since the eigenvalues $\gs_p(G_1(t_j))=\set{\pm i\hat{\gw}_k(t_j)}_{k=1}^q\cup \set{0}$ are distinct, it is easy to use the structures of $G_1(t_j)$ and the definition of $B_{1j}$ to show that $B_{1j}^\ast\phi\neq 0$ 
whenever $0\neq \phi\in \ker(G_1(t_j)^\ast)$ or
 $0\neq \phi\in \ker(\pm i\hat{\gw}_k(t_j)-G_1(t_j)^\ast)$ for some $k\in\List{q}$.
Thus the pair $(G_1(t_j),B_{1j})$ is controllable.
If we define $\gd=[\gd_1,\ldots,\gd_q]^T\in\R^q$ by
\eq{
\gd_1 &= \hat{\gw}_1(t_j), \quad
\gd_k = \hat{\gw}_k(t_j)-\hat{\gw}_{k-1}(t_j) \quad \forall k\in \List[2]{q},
}
we can define continuous functions $\tilde{G}_1:\R^q\to \Lin(\C^{p(2q+1)})$ and $\tilde{B}_1:\R^q\to \C^{p(2q+1)\times p}$ so that $\tilde{G}_1(\gd)=G_1(t_j)$ and $\tilde{B}_1(\gd)=B_{1j}$. 
The assumption that $(\hat{\gw}_k(t_j))_{k=1}^q$ satisfy $\EAC(M_f,\eps_f)$ implies that $\gd$ is contained in a compact set $\Omega\subset \R^q$ (determined by $\eps_f$, $M_f$, and $(A,B,C,D)$), and similarly as above, $(\tilde{G}_1(\gd),\tilde{B}_1(\gd))$ is controllable whenever $\gd\in\Omega$.
Because of this, Lemma~\ref{lem:IMstabRiccbound} implies that 
 there exist $M_G,M_{K1}>0$ (independent of $t_j$) such that
$\norm{K_1^j}\leq M_{K1}$ and $\norm{e^{ (G_1(t_j)+B_{1j}K_1^j)t}}\leq M_G e^{-rt}$ for all $t\geq 0$.
Morever, the
 definition $K^j=[K_1^j,K_{21}+K_1^jH_j]$ and $\norm{H_j}\leq M_H$ imply that $\norm{K^j}\leq M_K$ for some $M_K>0$ independent of $t_j$.

If  $\max_k\abs{\hat{\gw}_k(t_j)- \gw_k}\to 0$ as $j\to\infty$,
 then $(\gw_k)_{k=1}^\infty$ also satisfy $\EAC(M_f,\eps_f)$ and thus $K_1^\infty$ and $K^\infty$ are well-defined.
Clearly $G_1(t_j)\to G_1^\infty\in \Lin(Z_0)$ as $j\to\infty$ where $G_1^\infty$ is obtained by replacing $(\hat{\gw}_k(t_j))_k$ by $(\gw_k)_k$ in $G_1(t_j)$. Since $\gl\mapsto P_K(\gl)$ and 
$\gl\mapsto (\CL+DK_{21})R(\gl,A+BK_{21})$
 are continuous on $i\R$, also $B_{1j}\to B_1^\infty\in \C^{p(2q+1)\times m}$ and $H^j\to H^\infty \in \Lin(X,Z_0)$ as $j\to\infty$. We have 
from~\citel{Sun98}{Thm.~3.1}
 that $\Pi_{1j}$ depends continuously on $G_1(t_j)$ and $B_{1j}$, and therefore $\Pi_{1j}\to \Pi_1^\infty\in \Lin(Z_0)$ as $j\to\infty$. Thus the definitions of $K_1^j$ and $K^j$ imply that  $\norm{K^j-K^\infty}\to 0$ as $j\to\infty$.

Since $(A,B,C,D)$ is regular, the operator $H_j\in \Lin(X,Z_0)$ extends to $\ran(B)\subset X_{-1}$.
It is straightforward to check that $G_1(t_j) H_j=H_jA_K+G_2(\CL+DK_{21})$ and $B_{1j}=H_jB+G_2D$.
The definition $K^j=[K_1^j,K_{21}+K_1^jH_j]$ and
similar computations as in~\citel{PauPoh10}{Thm.~13} and~\citel{Pau16a}{Thm.~15} show 
\eq{
\MoveEqLeft\pmat{I&H_j\\0&-I}(A_s(t_j)+B_sK^j)\pmat{I&H_j\\0&-I}\\
&= \pmat{G_1(t_j)+B_{1j}K_1^j&0\\-BK_1^j&A+BK_{21}}
=: \tilde{A}_{sK}^j
}
with domain 
$\Dom(\tilde{A}_{sK}^j)=\setm{[z_1,z_2]^T\in Z_0\times X_B}{Az_2+B(K_{21}z_2-K_1^jz_1)\in X}$
(the domains of the block operators can be analysed as in the proof of~\citel{Pau16a}{Thm.~15}).
Fix $0<\ga_s<\min \set{r,-\gw_0(T_K(t))}$.
Denoting $R_G^j(\gl)=R(\gl,G_1(t_j)+B_{1j}K_1^j)$ and $R_K(\gl)=R(\gl,A+BK_{21})$,
 for all $\gl\in\C$ with $\re\gl>-\ga_s$ we have
\eq{
R(\gl,\tilde{A}_{sK}^j)
&= \pmat{R_G^j(\gl)&0\\-R_K(\gl)BK_1^jR_G^j(\gl)&R_K(\gl)}.
}
Since $\norm{R_G^j(\gl)}\leq M_G/(\re\gl +r)\leq M_G/(r-\ga_s)$ and since
$B$ is admissible with respect to the stable semigroup generated by $A+BK_{21}$, also
 $\sup_{\re\gl> -\ga_s}\norm{R_K(\gl)B}<\infty $~\citel{TucWei09book}{Prop.~4.4.6}.
Because of this, the Gearhart--Pr\"uss--Greiner theorem~\citel{EngNag00book}{Thm.~V.1.11} imply that 
there exists $\tilde{M}_s>0$ independent of $t_j$ such that
 $\norm{\tilde{T}_s^j(t)}\leq \tilde{M}_s e^{-\ga_s t}$ for all $t\geq 0$
(this uniform bound can be deduced, e.g., by applying~\citel{EngNag00book}{Thm.~V.1.11} to the semigroup $\diag(\tilde{T}_s^0(t),\tilde{T}_s^1(t),\ldots)$ on $\lp[2](Z_0\times X)$).
 Due to $\norm{H_j}\leq M_H$ and
 the similarity between $\tilde{T}_s^j(t)$ and $T_s^j(t)$ we finally have that there exists $M_s>0$ such that
$\norm{T_s^j(t)}\leq M_s e^{-\ga_s t}$ for all $t\geq 0$.

It remains to prove the existence of $M_B>0$. We have
\eq{
 R(\gl,A_s(t_j)+B_sK^j)B_s
&= 
\pmat{I&H_j\\0&-I}R(\gl,\tilde{A}_{sK}^j)\pmat{I&H_j\\0&-I}\pmat{G_2D\\B}\\
&= \pmat{I&H_j\\0&-I}\pmat{R_G^j(\gl)B_{1j}\\-R_K(\gl)BK_1^jR_G^j(\gl)B_{1j}-R_K(\gl)B}
}
for $\gl\in\C_+$.
Since $\norm{R_G^j(\gl)}\leq  M_G/r$, $\sup_{\gl\in\C_+}\norm{R_K(\gl)B}<\infty $, $\norm{H_j}\leq M_H$, $\norm{K_1^j}\leq M_{K1}$ and $\norm{B_{1j}}\leq M_{B1}$,
we indeed have $\norm{R(\gl,A_s(t_j)+B_sK^j)B_s}\leq M_B$ for all $\gl\in\C_+$ and for some $M_B>0$ independent of $t_j$.
\end{proof}


\end{document}